%% file: main.tex
\documentclass[final]{siamltex}
\include{config}

\title{Model order reduction for parametrized nonlinear hyperbolic problems as an application to Uncertainty Quantification}
\author{  R. Crisovan,  D. Torlo, R. Abgrall and S. Tokareva  \\
	Institut f\"ur Mathematik, Winterthurstrasse 190, \\CH 8057  Z\"urich, Switzerland 
}
\begin{document}
	\maketitle
\begin{abstract}	
	In this work, we focus on reduced order modeling (ROM) techniques for hyperbolic conservation laws with application in uncertainty quantification (UQ) and in conjunction with the well-known Monte Carlo sampling method. Because we are interested in model order reduction (MOR) techniques for unsteady non-linear hyperbolic systems of conservation laws, which involve moving waves and discontinuities, we explore the parameter-time framework and in the same time we deal with nonlinearities using a POD-EIM-Greedy algorithm \cite{Drohmann2012}. We provide under some hypothesis an error indicator, which is also an error upper bound for the difference between the high fidelity solution and the reduced one.
\end{abstract}



\pagestyle{myheadings}
\thispagestyle{plain}
\markboth{ R. Crisovan,  D. Torlo, R. Abgrall and S. Tokareva}{Model order reduction for parametrized nonlinear hyperbolic problems as an application to UQ}
\input{main_UQ}

\input{Bibliography}

\end{document}

%% file: config.tex
\usepackage[includehead, includefoot,hmargin={3.4cm,3.4cm},vmargin={2.85cm, 2.8cm}]{geometry}
\usepackage[latin1]{inputenc}
\usepackage{csquotes}
\usepackage[T1]{fontenc}
\usepackage{amsmath}
\usepackage{amsfonts}
\usepackage{subfigure}
\usepackage{subfigmat}
\usepackage{algorithm,algorithmic}
\usepackage{amssymb} 
\usepackage{url}
\newcommand{\R}{\mathbb R}
\newcommand{\Ppol}{\mathbb P}

\usepackage{graphicx}

\newcommand{\CC}{{\mathcal{ C}}}
\newcommand{\DD}{{\mathcal{ D}}}

\newcommand{\II}{{\mathcal{ I}}}

\newcommand{\LL}{{\mathcal{ L}}}

\newcommand{\OO}{{\mathcal{ O}}}
\newcommand{\PP}{{\mathcal{ P}}}
\newcommand{\QQ}{{\mathcal{ Q}}}

\newcommand{\UU}{{\mathcal{ U}}}

\newcommand{\WW}{{\mathcal{ W}}}

\def\varphibold{\boldsymbol{\varphi}}

\usepackage[usenames,dvipsnames]{color}
\usepackage{bm}

\usepackage[justification=centering]{caption}

\def\beq{\begin{equation}}
\def\eeq{\end{equation}}
\def\esplit{\end{split}}
\def\beqalign{\begin{array}{rl}}
\def\eeqalign{\end{array}}
\DeclareMathOperator*{\argmin}{\arg\!\min}
\DeclareMathOperator*{\argmax}{\arg\!\max}
\makeatletter
\newcommand{\ALOOP}[1]{\ALC@it\algorithmicloop\ #1%
  \begin{ALC@loop}}
\newcommand{\ENDALOOP}{\end{ALC@loop}\ALC@it\algorithmicendloop}

\makeatother

\def\EIM{\textrm{EIM}}
\def\POD{\textrm{POD}}
\def\RB{\textrm{RB}}
\def\dimEIM{N_\EIM}
\def\dimRB{N}
\def\truthN{N_h}

\def\dimtrain{N_{\text{train}}}

\newcommand{\uRB}[1]{\ubold_{\dimRB}^{#1}}
\newcommand{\utruth}[1]{\ubold_{h}^{#1}}

\def\Bbold{\mathbf{B}}

\def\fbold{\mathbf{f}}
\def\gbold{\mathbf{g}}
\def\hbold{\mathbf{h}}

\def\nbold{\mathbf{n}}

\def\qbold{\mathbf{q}}
\def\rbold{\mathbf{r}}

\def\ubold{\mathbf{u}}
\def\vbold{\mathbf{v}}
\def\wbold{\mathbf{w}}
\def\xbold{\mathbf{x}}

\def\deltabold{\boldsymbol{\delta}}
\def\mubold{\boldsymbol{\mu}}
\def\taubold{\boldsymbol{\tau}}

\def\varphibold{\boldsymbol{\varphi}}

\def\sigmabold{\mathbf{\sigma}}
\def\Phibold{\mathbf{\Phi}}

\def\Sigmabold{\mathbf{\Sigma}}

\def\0bold{\boldsymbol{0}}
\def\0{\mathbf{\0}}

\newcommand{\IE}{{\mathbb{E}}}
\newcommand{\IV}{{\mathbb{V}}}

\newcommand{\IN}{{\mathbb{N}}}

%% file: main_UQ.tex
\section{Introduction}\label{ch:intro_UQ}
Parametrized partial differential equations (PPDE) have received in the last decades an increasing amount of attention from research fields as engineering and applied sciences.  All these domains have in common the dependency of the PPDE on the input parameters, which are used to describe possible variations in the solution, initial conditions, source terms and boundary conditions, to name just a few. Hence, the solutions of these problems are depending on a large number of different input values, as in optimization, control, design, uncertainty quantification, real time query and other applications. In all these cases, the aim is to be able to evaluate in an accurate and efficient way an output of interest when the input parameters are varying. This will be very time consuming or can even become prohibitive when using high-fidelity approximation techniques, such as finite element (FE), finite volume (FV) or spectral methods. For this kind of problems, model order reduction (MOR) techniques are used, in order to replace the high-fidelity problem by one featuring a much lower numerical complexity. A key ingredient of MOR  are the reduced basis (RB) methods, which allow to produce fast reduced surrogates of the original problem by only combining a few high-fidelity solutions (\textit{snapshots}) computed for a small set of parameter values \cite{ito98,Rozza06reducedbasis,grepl2007}. The most common and efficient strategies available to build a reduced basis space are the proper orthogonal decomposition (POD) and the greedy algorithm. These two sampling techniques have the same objective but in very different approach forms: the POD method is most often applied only in one dimensional (1D) space and mostly in conjunction with (Petrov-)Galerkin projection methods, in order to build reduced-order models (ROM) of time-dependent problems \cite{POD,rathinam2003}, but also in the context of parametrized systems \cite{bui-thanh2004,bui-thanh2008,MartinKahlbacher2007,Tonn2010}. The disadvantage of this method is that it relies on the singular value decomposition (SVD) of a large number of snapshots, which might entail a severe computational cost. On the other side, greedy algorithm \cite{prudhomme2001,prud2002,rozza08} represents an efficient alternative to POD and is directly applicable in the multi-dimensional parameter domain. The algorithm is based on an iterative sampling from the parameter space fulfilling at each step a suitable optimality criterion that relies on a posteriori error estimates. 

A first challenge in the context of ROM deal with unsteady problems, so implicitly the exploration of a parameter-time framework is needed. In this case, the sampling strategy to construct reduced basis spaces for the time-dependent problem is POD-greedy \cite{haasdonk_pod_greedy2008} and is based on combining the POD algorithm in time, with a greedy algorithm in the parameter space. In general, all these methods are well suited for parametrized elliptic and parabolic partial differential equation models, for which their solutions are smooth with respect to the change of the inputs. We are interested instead, in parametrized hyperbolic systems of conservation laws, which involve moving waves and discontinuities such as shocks. It is well known that, in this case, the discontinuities will persist also in the parameter space, hence accurate surrogates have to be developed, in order to be able to capture the evolution of the discontinuous solutions. A second challenge refers to the nonlinear problems. In general, the computational efficiency of the RB method rely on affine assumptions, which is not the case for a big range of problems, including the hyperbolic ones.  Hence, in order to approximate nonaffine or nonlinear terms, one can make use of the empirical interpolation method (EIM) which approximates a general parametrized function by a sum of affine terms. This method was first introduced in \cite{barrault04} and in the context of ROM in \cite{grepl2007}. Some applications of the EIM method are discussed in \cite{Maday2009383} and an a posteriori error analysis is presented in \cite{grepl2007,EFTANG2010575}. There are only a few papers in the literature which are focused on MOR methods for parametric nonlinear hyperbolic conservation laws and they are based on: POD and Galerkin projection \cite{Rowley2004,kalashnikov2011}, domain partitioning \cite{taddei2015}, Gauss\text{-}Newton with approximated tensors (GNAT) \cite{gnat2013}, $L^1$-norm minimization \cite{AbgrallCrisovan2016,AbgrallCrisovan2018} or suitable algorithms extended to linear and nonlinear hyperbolic problems \cite{haasdonk_pod_greedy2008,Haasdonk2009}. The work of Drohmann, Haasdonk and Ohlberger \cite{Drohmann2012}, presents a new approach of treating nonlinear operators in the reduced basis approximations of parametrized evolution equations based on empirical interpolation namely, the PODEI-Greedy algorithm, which constructs the reduced basis spaces for the empirical interpolation in a synchronized way. 

In this paper, we focus on reduced order models for hyperbolic conservation laws based on explicit finite volume (FV) schemes. The FV schemes will be formulated within the framework of residual distribution (RD) schemes. The advantages of this alternative are: a better accuracy, a much more compact stencil, easy parallelization, explicit scheme and no need of a sparse mass matrix "inversion". For more details on RD, we refer to the work of Abgrall \cite{abgrall2006residual,abgrall_2012,abgrall2017}. However, we want to emphasize that our approach can be applied to any general FV formulation and RD is just our choice, and we have made this choice because we have an available code for free, so to speak. In this work, we concentrate on uncertainty quantification (UQ) applications for hyperbolic conservation laws. In practice, the input parameters are obtained by measurements (observations) and these measurements are not always very precise, involving some degree of uncertainty \cite{UQbook,UQbook_schwab}. A good example of hyperbolic conservation laws is when computing the flow past an airfoil or a wing, the inputs for this calculation, such as the inflow Mach number, the angle of attack, as well as the parameters that specify the airfoil geometry, are all measured with some uncertainty. This uncertainty in the inputs results in the propagation of uncertainty in the solution \cite{abgrall_mishra_uq}. Moreover, the need of model order reduction for UQ is obvious by just taking into account that these problems feature high-dimensionality, low regularity and arbitrary probability measures. However, the classical methods (Monte Carlo, stochastic Galerkin projection method, stochastic collocation method, etc) can not be applied directly to solve the underlying deterministic PDEs, since they might need millions of full solutions (or even more), in order to achieve a certain accuracy. Hence, with the help of reduced basis method, together with an a posteriori error estimate, we will be able to break the curse of dimensionality of solving high dimensional UQ problems whenever the quantities of interest reside in a low dimensional space. Up to our knowledge, there is no work done on hyperbolic conservation laws with applications in UQ and the only results that are available in literature are holding for elliptic PDEs \cite{chen_siam_2013,Chen2014,chen_uq_rb}.

In the first section we will present the problem of interest namely, the unsteady hyperbolic conservation laws and we will explain the RD scheme in relation with the nonlinear fluxes. In Section \ref{ch:algorithm} we will describe the algorithms that we are using for the construction of the reduced basis: POD-Greedy, PODEI. In Section \ref{ch:UQ} we describe the UQ method and in the last Section we present our numerical results.  

\section{Problem of interest}\label{ch:hyperbolic}

\subsection{Hyperbolic conservation laws}\label{sec:HypCons}
In this work, we consider high-dimensional models (HDM) arising from the space discretization of hyperbolic PPDEs. These problems are characterized by a parameter $\mubold \in \mathcal{P}$ from some set of possible parameters $\mathcal{P}\subset\R^{p}$. The unsteady problem then consists of determining the state variable solution $\ubold(\xbold,t;\mubold)$ on a bounded interval $D \subset \R^{d}, d=1,2,3$ and finite time interval $\R_{+}=[0,T],T>0$ such that the following system of $m, m\geq 1$ balance laws to be satisfied:
\begin{equation}\label{eq:hyp_cons}
\displaystyle{\left\{ 
	\begin{array}{l l }
	\displaystyle{ \ubold_{t}(\xbold,t;\mubold)} + \mathcal{L}(\xbold,t;\mubold)[\ubold(\xbold,t;\mubold)] &= \hbold(t;\mubold),~\xbold\in D,~t\in\R_{+}, \\
	\Bbold(\ubold;\mubold) &= \gbold(t;\mubold),~\xbold\in\partial D,~t\in\R_{+},\\
	\ubold(\xbold,t=0;\mubold) &= \ubold_0(\xbold;\mubold),~\xbold\in D,
	\end{array} \right.}
\end{equation}
where the operator $\mathcal{L}(\cdot,t;\mubold)=\text{div}f(\ubold(\xbold,t;\mubold))$ represents the divergence of the nonlinear flux $\fbold:\R^m\rightarrow(\R^m)^d$, $\Bbold$ is a suitable boundary operator, and $\hbold,\gbold$ are volume, respectively surface forces. Obviously, the moving shocks and discontinuities will depend on the different parameter settings $\mubold \in \mathcal{P}$ and will develop during time. The task of the RB method will be to capture the evolution of both smooth and discontinuous solutions.

The discrete evolution schemes are based on approximating high-dimensional discrete space $\mathcal{W}_h\subset L^2(D)$ (or subset of some Hilbert space), $\text{dim}(\mathcal{W}_h)=N_h$, where $h$ represents the characteristic mesh size and by approximating the exact solution at time-instances $0=t^0<t^1<\dots t^{K}=T$ i.e providing a sequence of functions $\ubold_h^k(\mubold):\R^{N_h}\rightarrow\R^m$ for $k=0,\dots,K$ such that $\ubold_{h}^{k}(\mubold)\approx\ubold(t_k;\mubold)$.

\subsection{Residual distribution scheme}\label{sec:RDscheme}
\begin{figure}[!htb]
	\minipage{0.55\textwidth}
	\includegraphics[width=\linewidth]{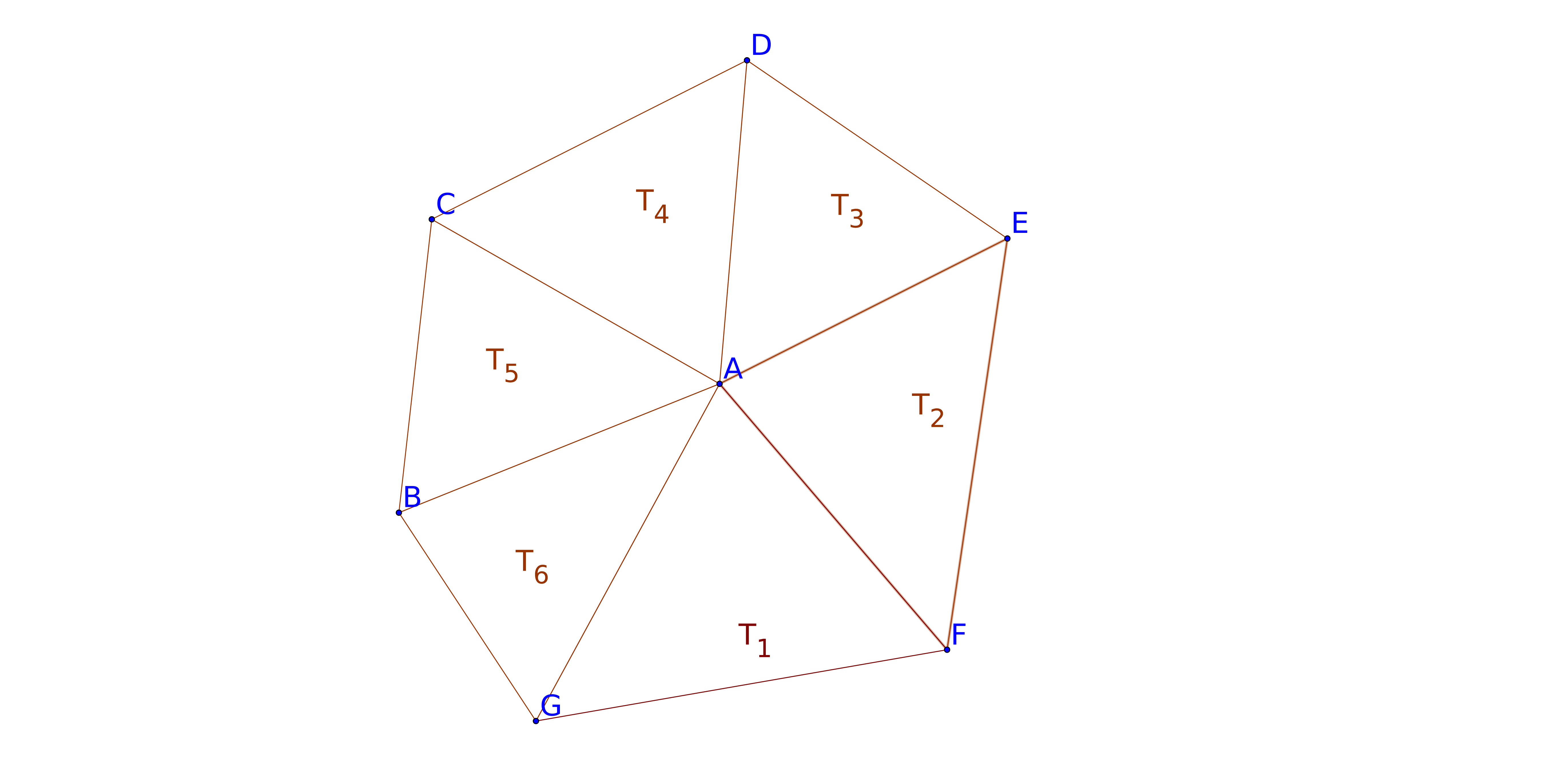}
	\caption{Triangulation $D_h$}\label{fig:triangulation}
	\endminipage\hfill
	\minipage{0.55\textwidth}
	\includegraphics[width=\linewidth]{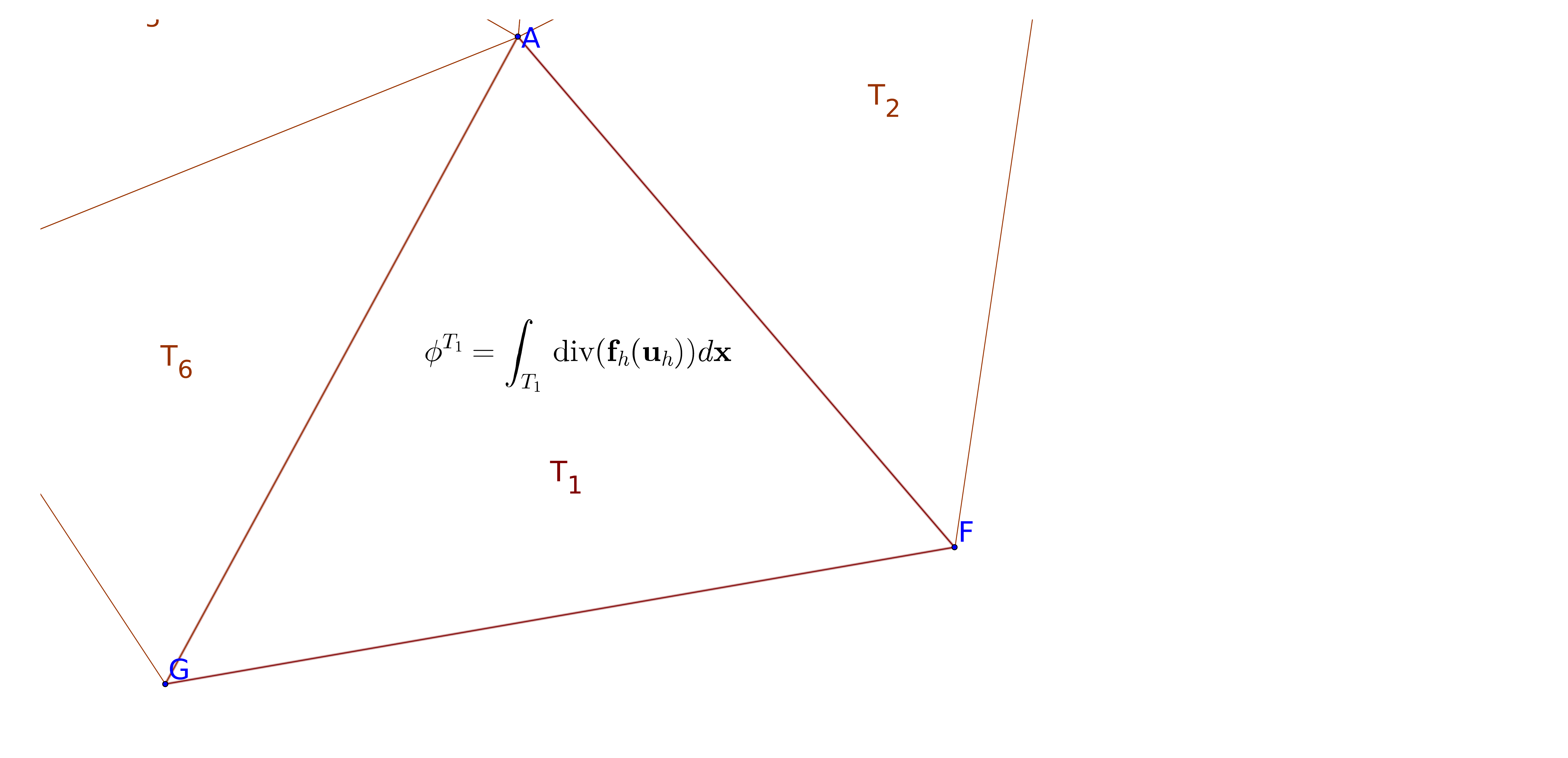}
	\caption{Compute the total residual}\label{fig:RD_fluctuation}
	\endminipage\hfill
	\minipage{0.55\textwidth}
	\includegraphics[width=\linewidth]{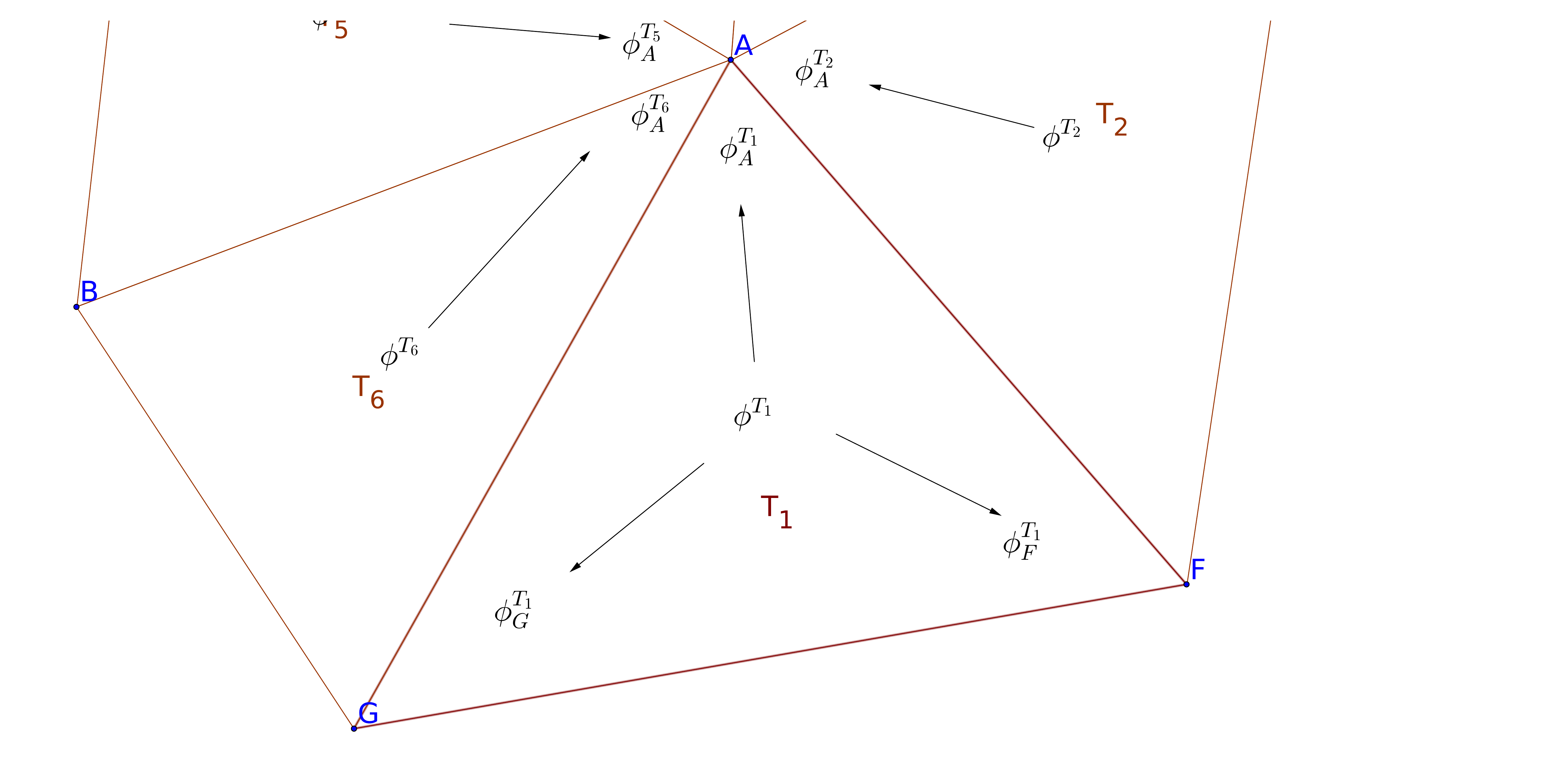}
	\caption{Compute the nodal residuals}\label{fig:RD_split_distribution}
	\endminipage\hfill
	\minipage{0.55\textwidth}%
	\includegraphics[width=\linewidth]{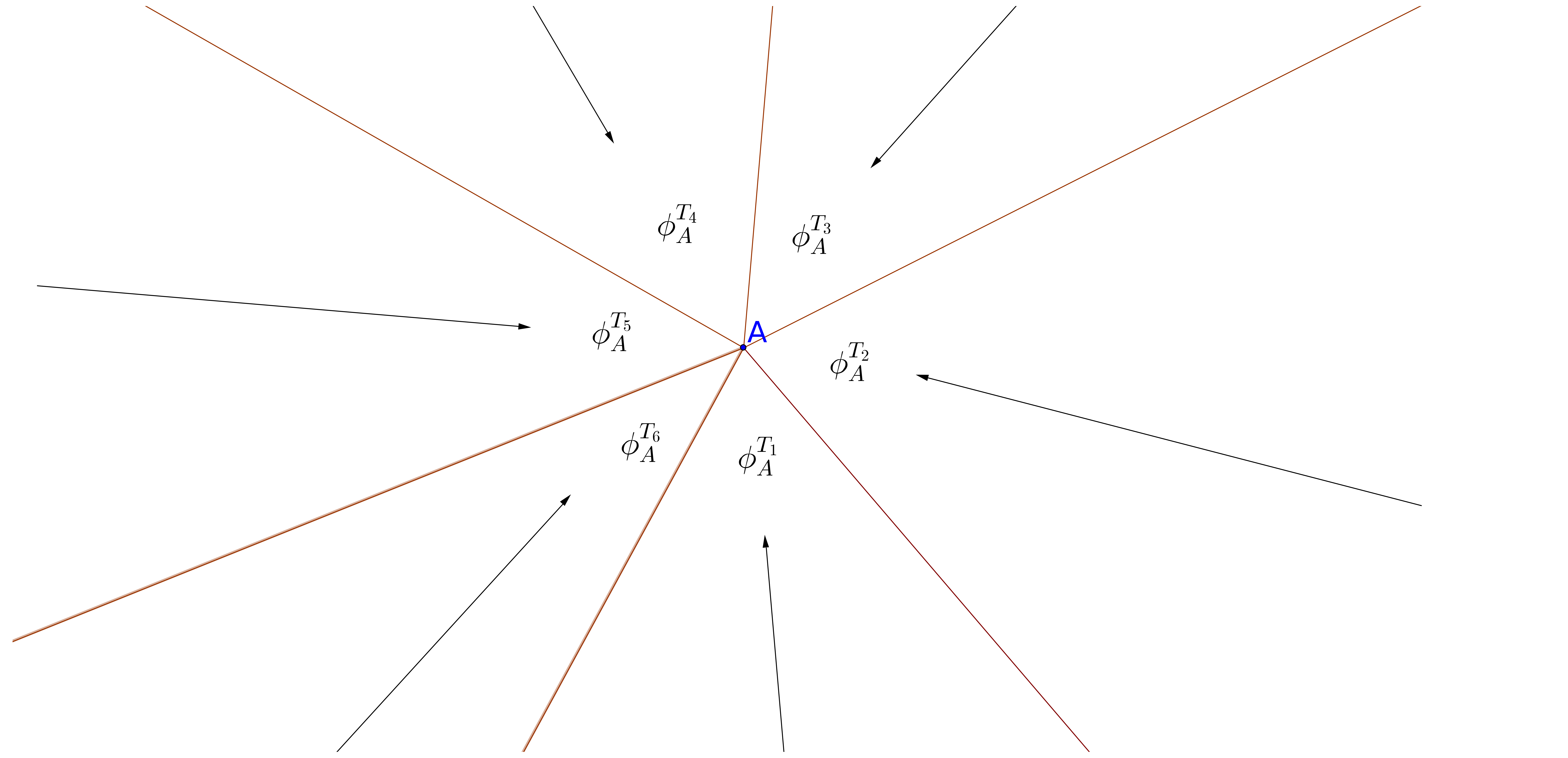}
	\caption{Collect all the residual contributions}\label{fig:RD_collect_contrib}
	\endminipage
\end{figure}

In this section, we are interested in the class of RD methods and we will show how any FV scheme can be written in this framework. We consider $D_h$ to be the triangulation of the domain $D$ (see Figure \ref{fig:triangulation}), $\Delta t_k=t_{k+1}-t_{k}$ the time steps for $k=0,\dots,K$ and we denote by $T$ a generic element of the mesh. We define the set $\sum_h:=\{\taubold_i\}_{i=1}^{N_h}\subset\mathcal{W}_{h}^{'}$ of linearly independent functionals, which are unisolvent on $\mathcal{W}_h$ i.e, there exist unique functions $\rho_i \in \mathcal{W}_h, ~i=1,\dots,N_h$ and satisfy:
$$\taubold_j(\rho_i)=\delta_{ij}, ~ 1\leq j \leq N_h.$$ The linear functionals $\taubold_i, i=1,\dots,N_h$ are called the degrees of freedom (DoFs) of the discrete function space $\mathcal{W}_h$, equipped with a scalar product $\langle \cdot, \cdot \rangle _{\WW_h}$ and a norm $||\cdot||_{\WW_h}$, and the functions $\rho_i,i=1,\dots,N_h$ are called the basis or shape functions. This shape functions can be for e.g, finite element, finite volume or discontinuous Galerkin basis functions on a numerical grid $D_h\subset D$.

In this case, the solution approximation space $\mathcal{W}_h$ is given by globally continuous piecewise polynomials of degree $r$:
\begin{equation}\label{eq:app_space}
\mathcal{W}_h=\{\ubold \in L^2(D_h)\cap C^{0}(D_h),\ubold_{|T}\in \Ppol^{r},\forall ~T \in D_h \}
\end{equation}
so that the numerical solution $\ubold_{h}^{k}$ can be written as a linear combination of shape functions $\rho_{i} \in \mathcal{W}_h, ~i=1,\dots,N_h$.

The main steps of the RD methods can be summarized as follows:
\begin{itemize}

	\item[1.] For any element $ T \in D_h$, compute the total residual
	\begin{equation}\label{TotalResidual_UQ}
	\Phi^{T}=\int_{T} \text{div }(\fbold_h(\ubold_h)) d\xbold=\int_{\partial T} \fbold_h(\ubold_{h}) \cdot \vec{\nbold}\ d\tilde{\xbold},
	\end{equation}
	where $\fbold_h$ is an approximation of $\fbold$ (Figure \ref{fig:RD_fluctuation}). 
	\item[2.] For any DoF $\taubold$ within an element $T$, define the nodal residuals $\Phi^{T}_{\taubold}$ as the contribution to the fluctuation term $\Phi^T$ (Figure \ref{fig:RD_split_distribution}) such that:
	\begin{equation}\label{DistributeTR_UQ}
	\sum_{\taubold \in T}\Phi_{\taubold}^{T}=\Phi^{T}.
	\end{equation}
	Equivalently, denoting by $\beta_{\taubold}^{T}$ the distribution coefficient of the DoF $\taubold$, we obtain:
	\begin{equation}\label{eq:distribution}
	\beta_{\taubold}^{T}=\frac{\Phi_{\taubold}^{T}}{\Phi^{T}}
	\end{equation}
	with
	\begin{equation}
	\sum_{\taubold \in T} \beta_{\taubold}^{T}=1.
	\end{equation}
	
	\item[3.] Assemble all the residual contributions $\Phi_{\taubold}^T$ from all elements $T$ surrounding a node $\taubold \in D_h$ (Figure \ref{fig:RD_collect_contrib}):
	\begin{equation}\label{AssembleTR}
	\sum_{T|\taubold \in T}\Phi_{\taubold}^{T}=0, ~\forall \taubold \in D_h.
	\end{equation}
\end{itemize}

This is a very general formulation and many classical schemes can be formulated within this framework. This variability hides mostly in how the residual of each triangle is distributed among the DoFs $\taubold \in T$, that is, on the choice of $\beta_{\taubold}^T$. For instance, distributing it evenly among nodes corresponds to a Lax-Friedriech type of scheme and can be defined without any reference to the geometry of a control volume, only by using the physical structure of the local flow. Another example, is the finite volume schemes, which are constructed using directions that are only related to the mesh definition and not to the structure of the solution. In this case, and whatever the order of accuracy of the scheme is, the approximation $\fbold_h(\ubold_h)$ is defined as the Lagrange interpolant of $\fbold(\ubold)$ at the DoF $\taubold \in T$.

\section{Algorithm}\label{ch:algorithm}
Before starting discussing the full algorithm we have used for our method, we should point out which are the main difficulties that we will encounter preparing our reduced basis space RB.\\
First of all, we know that the main prerequisite of a RB method is the separability into an affine decomposition, where the parameter dependent functionals are evaluated separately with respect to some precomputed parameter independent operators. To efficiently apply this principle to a non--linear functional, like our $\LL(\xbold,t;\mubold)[\ubold(\xbold, t; \mubold)]$, we need to introduce the empirical interpolation method in order to approximate an (\textit{a priori}) nonlinear parametrized operator with a separable one, which is efficient for evaluations of these operators for a reduced basis algorithm. We will show that this kind of surrogate operator can be computed in an efficient way using RD (or any FV) scheme in Section \ref{sec:online_phase}. Moreover, we need to build an efficient algorithm that will select sequentially some snapshots from some high--fidelity discretized solutions, until a prescribed tolerance. To do this, we will recur to a POD--Greedy algorithm, which is a combination of POD algorithm in time and a Greedy algorithm in the parameter space.\\
We will discuss in a general way the Greedy algorithm, since also EIM and POD--Greedy can be recast into a Greedy algorithm.
\subsection{Greedy algorithm}
A Greedy algorithm \cite{prudhomme2001,prud2002} is taking as an input some given precomputed functions and is building a reduced basis space, where the error of the approximation of any of these snapshots into this reduced basis space is smaller than a certain prescribed tolerance. The way the algorithm is choosing the reduced basis space, is an iterative method. At each step, the Greedy algorithm is selecting the snapshot that is worst approximated by the reduced basis projection and it is enriching the reduced basis space adding this new function. There are 3 main procedures that we will use in the Greedy algorithm. They are:
\begin{itemize}
	\item \textsc{InitBasis} which initializes the reduced basis $\DD_N$, also called dictionary in literature;
	\item \textsc{ErrorEstimate} which estimates the error between the high--fidelity function and its projection on the reduced basis space $\DD_N$;
	\item \textsc{UpdateBasis} which updates the RB space $\DD_N$, given a certain selected parameter.
\end{itemize}
The greedy algorithm proceeds as in Algorithm \ref{algo:greedy}.
\begin{algorithm}
	\fontsize{10pt}{10pt}\selectfont
	\caption{Greedy Algorithm} 
	\begin{algorithmic}[1]
		{\REQUIRE Training set $\mathcal{M}_{train} = \lbrace \mubold_i \rbrace_{i=1}^{\dimtrain}$, tolerance $\varepsilon^{tol}$ and $N_{max}$.
			\ENSURE Reduced basis $\DD_N$
			\STATE Initialize a reduced basis of dimension $N_0$:\\
			$\DD_{N_0}$= \textsc{InitBasis}\\
			$N=N_0$	
			\WHILE{$\max_{\mubold \in \mathcal{M}_{train} }  $\textsc{ErrorEstimate}$(\ubold(\mubold), \DD_N) \geq \varepsilon^{tol}$  AND $ N \leq  N_{max}$ } 
			\STATE Find the parameter of worst approximated snapshot:\\
			$\mubold_{max}=\argmax_{\mubold\in \mathcal{M}_{train}} \textsc{ErrorEstimate} (\ubold(\mubold), \DD_N)$
			\STATE Extend reduced basis $\DD_N$ with the found snapshot (adding the new snapshot to dictionary):\\
			$\DD_N, N=$\textsc{UpdateBasis}($\DD_N,\ubold (\mubold_{max})$)
			\ENDWHILE
		}
	\end{algorithmic}\label{algo:greedy}
\end{algorithm}
\subsection{Empirical Interpolation Method}\label{sec:EIM}
In this section we will apply the EIM algorithm \cite{barrault04} to the discretized operators. The method has the goal to apply an interpolation to the fluxes $\LL(\xbold,t^k;\mubold)[\ubold(\xbold, t^k; \mubold)] = \LL^k(\xbold, t^k ; \mubold)[\utruth{k}(\mubold)]$. The set of the interpolant DoFs $\Sigmabold_{\dimEIM}=\lbrace\taubold^{EIM}_m\rbrace_{m=1}^{\dimEIM}$, where $\taubold^{\EIM}_m \in \WW_h'$ and the corresponding set of interpolating basis functions $\QQ_{\dimEIM}=\lbrace\qbold_{m} \rbrace_{m=1}^{\dimEIM}$, where $\qbold_{m} \in \WW_h$ and $\taubold_m(\qbold_{n})=\deltabold_{mn}$ for $m\leq n$, will be the outputs of the algorithm. When the degrees of freedom can be identified with points in the domain (i.e. for Lagrange polynomial basis functions), EIM DoFs will be called ``magic points''. The specialization of Greedy algorithm into the EIM algorithm consists in the definition of the greedy procedures, i.e. Algorithm \ref{algo:EIM}, where the reduced basis, that we want to produce, comprise the interpolation DoFs $\Sigmabold_{\dimEIM}$ and the interpolation functions $\QQ_{\dimEIM}$, (i.e. $\DD_N=(\QQ_N, \Sigmabold_{N})$). After the EIM procedure, we will use the interpolated fluxes instead of the high fidelity discretized ones. 
\begin{equation}\label{eq:EIM_interpolation}
\II_{\dimEIM}[\LL (\xbold, t^k ; \mubold)][v_h]=\sum_{m=1}^{\dimEIM} \taubold_m^{\EIM}\left( \LL (\xbold, t^k ; \mubold)[v_h]\right) \qbold_{m} \approx \LL (\xbold, t^k ; \mubold)[v_h].
\end{equation}

\begin{algorithm}
	\fontsize{10pt}{10pt}\selectfont
	\caption{Empirical Interpolation Method} 
	\textsc{EIM--InitBasis()}
	\begin{algorithmic}[1]
		{	\RETURN empty initial basis $\DD_{0} = \emptyset $
		}
	\end{algorithmic}
	\hrule
	\vspace{1mm}
	\textsc{EIM--ErrorEstimate}($(\QQ_{M},\Sigmabold_{M}),\mubold, t^k$ )
	\begin{algorithmic}[1]
		{	\STATE Compute the exact flux $\vbold_h=\LL (\xbold, t^k ; \mubold)][\utruth{k}(\mubold)]$\\
			\STATE Compute the interpolation coefficients $\sigmabold^M(\vbold_h):= (\sigma^M_j)_{j=1}^M\in \R^M$ by solving the linear system (upper triangular)
			\begin{equation}\label{eq:EIM_system}
			\sum\limits_{j=1}^M \sigma_j^M(\vbold_h)\taubold_i^{\EIM}[\qbold_j]=\taubold_i^{\EIM}[\vbold_h],\quad \forall i=1,\dots,M
			\end{equation}
			\RETURN approximation error $|| \vbold_h- \sum_{j=1}^{M}\sigma_j^M(\vbold_h)\qbold_j||_{\WW_h}$
		}
	\end{algorithmic}
	\hrule
	\vspace{1mm}
	\textsc{EIM--UpdateBasis} ($(\QQ_{M},\Sigmabold_{M}), \mubold_{max}, t^{k_{max}}$)
	\begin{algorithmic}[1]
		{	\STATE Compute the exact flux\\ $\vbold_h=\LL (\xbold, t^{k_{max}} ; \mubold_{max})][\utruth{k_{max}}(\mubold_{max})]$\\
			\STATE Compute the interpolation coefficients \\ 
			$\sigmabold^M(\vbold_h):= (\sigma^M_j)_{j=1}^M\in \R^M$ from \eqref{eq:EIM_system} \\
			\STATE Compute the residual between the truth flux and its interpolant \\ 
			$\rbold_M=  \vbold_h- \sum_{j=1}^{M}\sigma_j^M(\vbold_h)\qbold_j$\\
			\STATE Find the DoF that maximize the residual \\ 
			$\taubold_{M+1}^{\EIM} := \argmax_{\taubold \in \Sigmabold_h} |\taubold (\rbold_M)|$\\
			\STATE Normalize the correspondent basis function\\
			$\qbold_{M+1}:=\taubold_{M+1}^{\EIM}(\rbold_M)^{-1} \cdot \rbold_M$ 
			\RETURN updated basis $\DD_{M+1}:=((\qbold_m)_{m=1}^{M+1}, (\taubold_{m}^{\EIM})_{m=1}^{M+1}) $
		}
	\end{algorithmic}\label{algo:EIM}
\end{algorithm}

The algorithm produced a basis $\QQ_{\dimEIM}$ which fulfills in a relaxed way the Kronecker's delta condition: $\taubold_m^{\dimEIM}(\qbold_n) =\deltabold_{mn}$ only if $m\leq n$. This condition will provide an upper triangular matrix that can be easily inverted during the EIM procedure to solve the interpolant coefficients problem. Moreover, the EIM basis functions spaces will be hierarchical, i.e. $\QQ_M\subset \QQ_{M+1}$, and the infinity norm of all the basis functions will be equal to 1 ($||\qbold_{m}||_{\infty}=1$). \\ 
Let us remark that, when we are dealing with Lagrange polynomial basis functions, formula \eqref{eq:EIM_interpolation} requires the evaluation of functions $\LL(\xbold, t^k ; \mubold)[v_h]$ only in the \textit{magic points}, and this will give the biggest reduction in computational time, since the evaluation of fluxes can be very expensive. Indeed, the number of interpolation DoFs should be $\dimEIM \ll N_h$. In RD framework, we can explicitly see what we need to compute:
\begin{equation}
\taubold_i [\LL (\xbold, t^k ; \mubold)][\utruth{k}(\mubold)] = \sum\limits_{T|i \in T} \Phibold^T_i (\utruth{k}(\mubold)).
\end{equation} 
Each nodal residual $\Phibold^T_i (\utruth{k}(\mubold))$ depends only on DoFs of element T, this means that for each \textit{magic point} $i$ we have to keep track of the function $\utruth{k}(\mubold)$ in all the DoFs of the elements $T$  to which $i$ belongs. The number of these DoF is \textit{mesh--dependent}, for the simplest example in 1D with $\mathbb{P}^1$ piecewise continuous elements we know that for each magic point we have to keep track of 3 points: itself, its right and left neighborhoods. If we suppose some regularities on the mesh we can say that at most each vertex belongs to $C$ elements. In this case, again for $\mathbb{P}^1$ Lagrangian basis functions, the number of DoF we are interested in is $R=C (K-2) +1$, where $K$ is the biggest number of vertices that an element $T$ can have.\\
At the end, we will have that the empirical interpolation method will provide an approximated version of the fluxes that depends at most on $R\dimEIM \ll N_h$ DoFs. 
\subsection{POD--Greedy}\label{sec:pod_greedy}
To create a reduced basis $\RB$ space, we want to find a low dimensionality \textit{good} approximation of the high fidelity functional space $\WW_h$. The algorithm that will provide this is a combination of different algorithms, such as POD \cite{POD, PCA}, POD--greedy \cite{RB_book_rozza}, EIM--greedy \cite{barrault04}. What we will get is a POD--EIM--greedy algorithm, described by \cite{Drohmann2012}. The main idea is to extend EIM basis functions and POD--greedy basis functions in a synchronized way, at each step of the main greedy algorithm.

A key ingredient of the procedure is the POD method, which is also known as PCA (principal component analysis) in statistical environment. The POD receives as input a set of vectors and returns the subspace of dimension $N_{\POD}$ which best represents the vectors given as a projection onto this subspace. We can write it in this way
\begin{equation}\label{def:POD}
POD(\lbrace \vbold_i \rbrace_{i=1}^N) = \argmin\limits_{U| \text{dim}(U)=N_{\POD}} \max_{i \in \lbrace 1, \dots, N \rbrace} \left(||\vbold_i-\PP_U(\vbold_i) ||_2\right).
\end{equation}
Equivalently, this can be seen as the subspace of fixed dimension that maximizes the variance. The algorithm is based on SVD decomposition. We need to order the eigenvalues from the biggest to the smallest and we keep the first $N_{\POD}$ ones and the related eigenvectors. The span of the latter will be the output of the algorithm. To choose the dimension of this subspace, it is possible to use a tolerance, which will decide which percentage of the variance we want to keep or which percentage of the error we want to ignore. In our algorithms, we will use different tolerances, according to whether we want them to be fast (bigger $N_{\POD}$) or sharp (small $N_{\POD}$, even 1). 

Before explaining the main algorithm, let us introduce the POD-Greedy algorithm, which deals with unsteady problems in the reduced basis context. The goal of the algorithm is to select new basis functions iteratively between precomputed snapshots $\lbrace\lbrace \utruth{k}(\mubold_i) \rbrace_{k=1}^K \rbrace_{i =1}^{\dimtrain}$. So, we have to find strategies to go through the parameter space and through the time steps. First, we explore the parameter space through a Greedy algorithm. We pick the parameter $\mubold_{max}$ that is worst approximated in RB space. Hence, on its temporal evolution $\lbrace \utruth{k}(\mubold_{max})\rbrace_{k=1}^K$, we perform a POD that chooses the most representative $M$--dimensional space for that solution, to compress the solution in a few synthetic basis functions. Then we add to the RB space the new basis functions selected by POD. Finally, we perform a second POD on the RB space, to get rid of useless information.\\
Overall, we will compute a Greedy algorithm on the parameter domain $\mathcal{P}$ and a POD on the temporal space. Also in this case, we can write the POD-Greedy Algorithm \ref{algo:POD_greedy}, specifying the greedy procedures as in Algorithm \ref{algo:greedy}.
\begin{algorithm}
	\fontsize{10pt}{10pt}\selectfont
	\caption{POD--Greedy} 
	\textsc{POD--Greedy--InitBasis()}
	\begin{algorithmic}[1]
		{	\STATE Pick a parameter $\mubold$ and compute the solution through all the time steps $t^k$: $ \lbrace \utruth{k}(\mubold) \rbrace_{k=1}^K$\\
			\RETURN initial basis $\DD_{0} = \POD (\lbrace \utruth{k} (\mubold) \rbrace_{k=1}^K )$
		}
	\end{algorithmic}
	\hrule
	\vspace{1mm}
	\textsc{POD--Greedy--ErrorEstimate}($\RB, \mubold, t^k$ )
	\begin{algorithmic}[1]
		{	\RETURN error indicator $\eta^k_{\dimRB ,\dimEIM} (\mubold) \geq ||\utruth{k}(\mubold)-\uRB{k} (\mubold) ||_{\WW_h}$
		}
	\end{algorithmic}
	\hrule
	\vspace{1mm}
	\textsc{POD--Greedy--UpdateBasis} ($\RB, \mubold_{max} $)
	\begin{algorithmic}[1]
		{	\STATE Compute the exact solution for all timestep with high fidelity solver  $\lbrace \utruth{k}(\mubold_{max}) \rbrace_{k=1}^K  $\\
			\STATE Compute the Galerkin projection of the solution onto the $\RB$ space $\PP[\utruth{k}(\mubold_{max})]$
			\STATE Compute the POD over time steps of the orthogonal projection of the high fidelity solution \\ 
			$\RB_{add}=\POD\left( \lbrace \PP[\utruth{k}(\mubold_{max})]-\utruth{k}(\mubold_{max})\rbrace_{k=1}^K\right)$\\
			\STATE Compute a second POD to get rid of extra information \\ 
			$\RB=\POD(\RB_{add} \cup \RB)$\\
			\RETURN updated basis $\RB $
		}
	\end{algorithmic}\label{algo:POD_greedy}
\end{algorithm}

Let us point out a couple of details of Algorithm \ref{algo:POD_greedy}. At the beginning, we may initialize the reduced basis with a {\POD}  with a $N_{\POD}$ bigger than one used later (or a smaller error tolerance), since we still do not have any $\RB$ and we want to accelerate the first steps, to decrease the number of greedy steps. During the rest of the algorithm we will use the $\POD$ on the time evolution of the worst approximated solution in the training set and $N_{\POD}$ here will be smaller (or the tolerance will be bigger). The last $\POD$ that we use is in the last step of the \textsc{POD--Greedy--UpdateBasis}, where $N_{\POD}$ will be big and set by a very small tolerance (of the order of the final error that we want to reach). This will kill some spurious vectors that may come from oscillations or small errors. Often this step is not changing the updated reduced basis. 

About the error estimator $\eta$, we would like to have a function which is independent of $\truthN$ that can be computed also in an \textit{online phase}. Of course, this bound should also be enough sharp, to give a precise idea of the error. We will describe in section \ref{sec:error_indicator}, an error indicator that is possible to use. If this indicator is not available, in the \textit{offline phase} we can still use the real error, which is computationally less efficient, and in the \textit{online phase}, where the high fidelity solutions are not available, we can not compute it directly. So, we will not have an explicit error bound to guarantee a good approximation. 

In Algorithm \ref{algo:POD_greedy}, it is not written explicitly the EIM--method that every time we are applying to some reduced basis solutions. Moreover, the error indicator should also include the error produced by EIM procedure. This approach has some drawbacks described in \cite{Drohmann2012}: 
\begin{enumerate}
	\item Is not really clear what is the relation between the tolerance used to stop EIM algorithm and the error produced in the POD--Greedy and how it influences the error indicator $\eta$. Therefore, it is impossible to determine a priori an optimal correlation between the reduced basis space and the EIM space.
	\item The empirical interpolation error estimation depends on high dimensional computations for each parameter and time step tested. This can be very inefficient.
\end{enumerate}

\subsection{PODEIM--Greedy}
To avoid these drawbacks, the idea of \cite{Drohmann2012} is to synchronize the EIM and the POD--Greedy algorithms. We sketch the steps of the PODEIM-Greedy in Algorithm \ref{algo:PODEIM_greedy} with the remark that also this algorithm can be rewritten in terms of a greedy one \ref{algo:greedy}.

\begin{algorithm}
	\fontsize{10pt}{10pt}\selectfont
	\caption{PODEIM--Greedy} 
	\textsc{PODEIM--Greedy--InitBasis()}
	\begin{algorithmic}[1]
		{	\STATE $(\QQ_{M_{small}} , \Sigmabold_{M_{small}}) = \textsc{EIM-Greedy}(\mathcal{M}_{train}, \varepsilon_{tol,small})$
			\STATE Pick a parameter $\mubold$ and compute the solution through all the time steps $t^k$: $ \lbrace \utruth{k}(\mubold) \rbrace_{k=1}^K$\\
			\STATE $\RB_{0} = \POD (\lbrace \utruth{k} (\mubold) \rbrace_{k=1}^K )$
			\RETURN  initial bases $\DD_{0}=(\RB_{0} , (\QQ_{M_{small}} , \Sigmabold_{M_{small}}))$
		}
	\end{algorithmic}
	\hrule	
	\vspace{1mm}
	\textsc{PODEIM--Greedy--ErrorEstimate}($\DD_{S}, \mubold, t^k$ )
	\begin{algorithmic}[1]
		{	\RETURN error indicator $\eta^k_{\dimRB ,\dimEIM} (\mubold)$
		}
	\end{algorithmic}
	\hrule
	\vspace{1mm}
	\textsc{PODEIM--Greedy--UpdateBasis} ($\DD_{S}, \mubold_{max} $)
	\begin{algorithmic}[1]
		{	\STATE Extend EIM basis $D^{\EIM}_{\dimEIM+1} = $ \textsc{EIM--UpdateBasis} $(D^{\EIM}_{\dimEIM}, \mubold_{max})$
			\STATE Extend RB basis $D^{\RB}_{\dimRB+1} =$ \textsc{POD--Greedy--UpdateBasis} $(D^{\RB}_{\dimRB}, \mubold_{max})$ \\
			\STATE Discard extended {\RB} if error increases:
			\IF{ $\eta^k_{\dimRB-1,\dimEIM -1}(\mubold_{max})< \max_{\mubold_i \in \mathcal{M}_{train} } \eta^k_{\dimRB,\dimEIM }$} 
			\RETURN only {\EIM} updated basis: $\DD_{S+1} = (D^{\RB}_{\dimRB},D^{\EIM}_{\dimEIM+1})$\\
			\ELSE 	\RETURN updated basis $\DD_{S+1} = (D^{\RB}_{\dimRB+1},D^{\EIM}_{\dimEIM+1})$	\\
			\ENDIF
		}
	\end{algorithmic}\label{algo:PODEIM_greedy}
\end{algorithm}

The differences between this new algorithm and the POD--Greedy are in the update phase, where we enrich at the same moment the EIM and the {\RB} basis. Moreover, it is possible that the error (and the indicator $\eta$) is not monotonically decreasing as the dimension of {\RB} increases. This is caused by a bad approximation of the non--linear fluxes through the EIM. Indeed, in such a situation, we are enlarging only the EIM space and discarding the additional part of the {\RB} that we added. This leads to an automatic tuning between $\dimRB$ and $\dimEIM$. 

\subsection{Online--phase}\label{sec:online_phase}
In this section we will describe the reduced basis scheme that we will eventually apply to find a reduced solution. This process is also used in the \textit{offline--phase} at each greedy step for each parameter in the training set, to get the reduced solution and the correspondent error. We will focus on explicit finite volume method, that can be rewritten into RD explicit scheme, but it is possible to extend this scheme to implicit (Newton iteration based method) as done in \cite{Drohmann2012}. The basic idea is to replace the discrete evolution operator $\LL[\cdot]:=\LL(\xbold,t^k;\mubold)[\cdot]$ with its empirical interpolants and project it onto the $\RB$ space. For this purpose, let us introduce the orthogonal projection $\Pi :\WW_h \to \RB  $ such that
\begin{equation}
\langle \Pi [u], \varphibold \rangle_{\WW_h}=\langle u, \varphibold \rangle_{\WW_h} , \qquad \forall \varphibold \in \RB
\end{equation}
and we can define the reduced operator as 
\begin{equation}
\LL_{\RB}:=\Pi \circ \II_{\dimEIM} \circ \LL .
\end{equation}
Let us define $\lbrace \varphibold_{\RB,i} \rbrace_{i =1}^{\dimRB}$ a basis of $\RB$, $\lbrace\qbold_{m} \rbrace_{m=1}^{\dimEIM}$ the interpolation functions of {\EIM} space and, for $m=1,\dots,\dimEIM$, let us define $\lbrace \theta_i^m\rbrace_{i=1}^{\dimRB}$ such that $\Pi(\qbold_{m})=\sum_{i=1}^{\dimRB} \theta_i^m \varphibold_{\RB,i} $.

To begin the procedure, for any parameter $\mubold$, we compute the trajectory of the reduced solution, projecting the initial data onto the {\RB} space: $\uRB{0}(\mubold):=\Pi[\utruth{0}(\mubold)]$. Then, for each time step, we compute the reduced solution applying the reduced operator $\LL_{\RB}[\uRB{k}]$. This implies to compute 

\begin{equation}
\begin{split}
\uRB{k+1}(\mubold)= &\uRB{k}(\mubold)-\LL_{\RB}[\uRB{k}(\mubold)] = \sum_{i=1}^{\dimRB } \alpha_{\RB,i}^k(\mubold) \varphibold_{\RB,i} -\Pi( \II_{\dimEIM} ( \LL [\uRB{k}(\mubold)]))=\\
=& \sum_{i=1}^{\dimRB } \alpha_{\RB,i}^k(\mubold) \varphibold_{\RB,i}-\Pi\left( \sum_{m=1}^{\dimEIM} \taubold_m^{\dimEIM}( \LL [\uRB{k}(\mubold)])\qbold_{m}\right)=\\
=&\sum_{i=1}^{\dimRB } \alpha_{\RB,i}^k(\mubold) \varphibold_{\RB,i}- \sum_{m=1}^{\dimEIM} \taubold_m^{\dimEIM}( \LL [\uRB{k}(\mubold)]) \Pi \left(\qbold_{m}\right)=\\
=&\sum_{i=1}^{\dimRB } \alpha_{\RB,i}^k(\mubold) \varphibold_{\RB,i} - \sum_{m=1}^{\dimEIM} \taubold_m^{\dimEIM}( \LL [\uRB{k}(\mubold)]) \sum_{i=1}^{\dimRB } \theta_i^m \varphibold_{\RB,i} =\\
=&\sum_{i=1}^{\dimRB } \left(\alpha_{\RB,i}^k(\mubold) - \sum_{m=1}^{\dimEIM} \taubold_m^{\dimEIM}( \LL [\uRB{k}(\mubold)])  \theta_i^m \right)\varphibold_{\RB,i}.
\end{split}
\end{equation}

In the last formula, what we really need to compute \textit{online} is only $\taubold_{m}( \LL [\uRB{k}(\mubold)]) ,\, \forall m=1,\dots , \dimEIM$, which implies, as written in Section \ref{sec:EIM}, $R\dimEIM$ evaluation of the flux. All the other terms are computed previously and stored: $\alpha_{\RB,i}^k(\mubold)$ are the coefficient of the previous time step, $\varphibold_{\RB,i}$ are the basis functions of $\RB$, previously computed, and $\theta_i^m$ are the projection coefficient of {\EIM} functions onto $\RB$.  Overall, the computational cost of a reduced solution at each time step will be $\OO (R\dimEIM)$ flux evaluations and $ \OO(\dimEIM \dimRB)$ multiplications.

\subsection{Error indicator} \label{sec:error_indicator}
We can provide an error indicator, which is also an error upper bound for the difference between the high fidelity solution and the reduced one, under some hypothesis. This estimation is derived following the guidelines of \cite{Drohmann2012} and \cite{Haasdonk2009}. The hypothesis under which the indicator becomes a bound is that there exists a higher order empirical interpolation of the used operators which is exact. This requirement is fulfilled if we take the interpolation over all the DoFs ($\dimEIM': \dimEIM +\dimEIM'= H $), where $H$ is the number of DoFs. But, for practical purposes, it has been show in \cite{Drohmann2012} that fewer points are necessary to get a good estimator. \\
Let us define other $\dimEIM'$ {\EIM } basis functions $\lbrace \qbold'_{m} \rbrace_{m=1}^{\dimEIM'}$, simply iterating further the {\EIM} procedure. And we suppose that 
\begin{equation}\label{eq:EIM_M1_hypothesis}
\II_{\dimEIM+\dimEIM'}[\LL(\xbold, t^k ; \mubold)][\uRB{k}(\mubold_i)] = \LL(\xbold, t^k ; \mubold)[\uRB{k}(\mubold_i)].
\end{equation}
Moreover, we suppose that the projection of the initial condition are in the reduced basis space, i.e. $\utruth{0}(\mubold)\in\RB,\, \forall \mubold \in \PP$. This can be easily obtained if there exists an affine decomposition of the parametric dependent part of the initial conditions: $\utruth{0}(\xbold,\mubold)=\sum \limits_{k=1}^F \alpha_k(\mubold) u_k(\xbold)$. Anyway, we will show that, also without fulfilling this condition, the numerical results do not present particular problems if the tolerance of the $\RB$ is enough small.\\
Then, we need a very last hypothesis on the operator $\textrm{Id} - \Delta t \LL (\xbold, t^k ; \mubold)$ namely, to be Lipschitz continuous with constant $C>0$, i.e. $\forall u,v \in \WW_h$: 
\begin{equation}\label{eq:Lipschitz}
||u-v -\Delta t \LL [u] + \Delta t \LL [v]||_{\WW_h}  \leq C ||u-v||_{\WW_h}
\end{equation}
holds.\\
Under these hypothesis we can say that the error $e^k(\mubold) := \utruth{k}(\mubold) - \uRB{k}(\mubold)$ can be bounded by $\eta^k_{\dimRB,\dimEIM,\dimEIM'}(\mubold)$, which can be computed efficiently, and it is defined as 
\begin{equation}\label{eq:error_estimator}
||e^K(\mubold)||_{\WW_h} \leq \eta^K_{\dimRB,\dimEIM,\dimEIM'}(\mubold):= \sum_{k=1}^K C^{K-k} \left( \sum_{m=1}^{\dimEIM'} \Delta t \theta^k_m (\mubold)  \left\lVert  \qbold'_{m} \right\rVert_{\WW_h}+ \Delta t || R^k(\mubold)||_{\WW_h} \right),
\end{equation}
where 
\begin{equation}
\Delta t R^k(\mubold) := \uRB{k}(\mubold) -\uRB{k-1}(\mubold)+ \Delta t \II_{\dimEIM} [\LL] [\uRB{k-1}(\mubold)] \end{equation}
and the coefficient 
\begin{equation}
\theta^k_m(\mubold)= \taubold_m^{\dimEIM'} \left( \LL [\uRB{k-1}(\mubold)] \right), \, \forall m= 1, \dots, \dimEIM'.
\end{equation}
\begin{proof}
	For the sake of simplicity, we will drop all the $\mubold$ parameters.
	\begin{equation}
	\begin{split}
	\left\lVert \utruth{K+1}- \uRB{K+1} \right\rVert =& \left\lVert(\textrm{Id} - \Delta t \LL) (\utruth{K}) - (\textrm{Id} - \Delta t  \II_{\dimEIM} [\LL]) (\uRB{K} ) - \Delta t R^{K} \right\rVert =\\
	\leq & \left\lVert (\textrm{Id} - \Delta t \LL) (\utruth{K}) -  (\textrm{Id} - \Delta t \LL) (\uRB{K})\right\rVert+\left\lVert (\Delta t \LL-  \Delta t  \II_{\dimEIM} [\LL])(\uRB{K} )\right\rVert +  \\ & + || \Delta t R^{k} || . 
	\end{split}
	\end{equation}
	Then we can use  Lipschitz condition \eqref{eq:Lipschitz} and get the following:
	\begin{equation}
	\begin{split}	
	\left\lVert \utruth{K+1}- \uRB{K+1} \right\rVert \leq & C \left\lVert \utruth{K} -  \uRB{K}\right\rVert+\left\lVert (\Delta t \LL-  \Delta t  \II_{\dimEIM} [\LL])(\uRB{K} )  \right\rVert + || \Delta t R^{K} ||. 
	\end{split}
	\end{equation}
	Now, using the fact that the evolution is exactly represented with the second EIM interpolant \eqref{eq:EIM_M1_hypothesis}, we can rewrite it into:
	\begin{equation}
	\begin{split}	 
	C \left\lVert \utruth{K} -  \uRB{K}\right\rVert+ & \left\lVert (\Delta t \II_{\dimEIM+\dimEIM'}[\LL]-  \Delta t  \II_{\dimEIM} [\LL])(\uRB{K} )  \right\rVert + || \Delta t R^{K} ||\leq\\  
	\leq & C \left\lVert \utruth{K} -  \uRB{K}\right\rVert+\left\lVert \Delta t \sum_{m=1}^{\dimEIM'} \taubold_m^{\dimEIM'} [\LL(\uRB{K} )] \qbold'_{m}  \right\rVert + || \Delta t R^{K} ||\leq\\
	\leq & C \left\lVert \utruth{K} -  \uRB{K}\right\rVert+\left\lVert \Delta t \sum_{m=1}^{\dimEIM'} \theta^K_m \qbold'_{m}  \right\rVert + || \Delta t R^{K} ||\leq\\
	\leq & \sum_{k=1}^{K+1} C^{K+1-k} \left( \left\lVert \sum_{m=1}^{\dimEIM'} \Delta t \theta^k_m (\mubold)  \qbold'_{m} \right\rVert + ||\Delta t R^k(\mubold)|| \right).
	\end{split}
	\end{equation}
	This proves that the error indicator is an actual bound when all the hypothesis are fulfilled.
\end{proof}
\par Anyway, from experimental results, we can see that, also when we are not in this case, the estimator is giving a good approximation of the error. Indeed, for $\EIM'$, as shown in \cite{Drohmann2012}, we can take very few basis functions and get good results, because the chosen DoFs should be the ones that maximize the error. Moreover, its computational cost is $\OO (R\dimEIM')$ evaluations of the flux. 
\subsubsection*{Estimation of the Lipschitz constant}
A couple of words should be spent on the way to find the Lipschitz constant $C$. Actually, it really depends on the specific method that is used and it is difficult to give a general way to estimate it. For the scheme that we use, we could not find a sharp estimation, because it involves some operators that do not belong to $\CC^1$. But, since the operator $\LL$ is the discretized operator of the gradient of the flux, we can use the spectral radius $\rho$ of the Jacobian of the flux to approximate this constant. 
\begin{equation}
\begin{split}
&||u-v - \LL[u] + \Delta t \LL [v]||\approx ||u-v||+\Delta t||\nabla f(u) - \nabla f(v)||  \approx \\
\approx &||u-v|| + \Delta t|| J(f) (u-v)||\leq ||u-v||+\rho \Delta t||u-v|| = (1+\rho \Delta t) ||u-v||.
\end{split}
\end{equation} 
What we used in the numerical experiments is a bound $b$ for the spectral radius of the Jacobian of the flux, for $u$ being in a reasonable box. Then we can fix $C=1+b \Delta t$. This can be done in a smarter way and more efficiently if the flux is affinely depending on the parameter $\mubold$. Therefore, one can split this constant into a parameter dependent and a fixed part.



\section{Applications to Uncertainty Quantification}\label{ch:UQ}

\subsection{Stochastic conservation laws}

Many problems in physics and engineering are modeled by hyperbolic
systems of conservation or balance laws.  As examples for these equations,
we can mention the Euler equations of compressible gas dynamics, the Shallow Water Equations of hydrology, the Magnetohydrodynamics
(MHD) equations of plasma physics, see, e.g.~\cite{God,Dafermos}. 

Many efficient numerical methods have been developed to approximate the entropy solutions
of systems of conservation laws \cite{God,LeVeque2002}, e.g. finite volume 
or discontinuous Galerkin methods. The classical assumption in designing efficient numerical
methods is that all the input data, e.g. initial and boundary conditions, flux vectors, sources, etc, are deterministic. However, in many situations of practical interest, these data are subject to inherent uncertainty in modeling and measurements of physical parameters.
Such incomplete information in the uncertain data can be represented mathematically as random fields.
Such data are described in terms of 
statistical quantities of interest like the mean, variance, 
higher statistical moments; in some cases the distribution 
law of the stochastic data is also assumed to be known. 

A mathematical framework of \textit{random entropy solutions} 
for scalar conservation laws with random initial data has been developed in \cite{MishraSchwab2010}.
There, existence and uniqueness of random entropy solutions
has been shown for scalar hyperbolic conservation laws,
also in multiple dimensions. 
Furthermore, the existence of the statistical quantities 
of the random entropy solution such as the statistical mean 
and $k$-point spatio-temporal correlation functions under 
suitable assumptions on the random initial data have been proven. The existence and uniqueness of the random entropy solutions for scalar conservation laws with random fluxes has been proven in \cite{SchwabMishraRisebroTokareva2012}.

A number of numerical methods for uncertainty quantification (UQ) in hyperbolic conservation laws 
have been proposed and studied recently in e.g.
\cite{MishraSchwab2010,MishraSchwabSukys2011,Troyen-1,Troyen-2,Abgrall,Poette,Lin-Su-Karniadakis-1,Lin-Su-Karniadakis-2,Gottlieb-Xiu,SchwabTokareva2011,SFVSpringer2014}.

\subsection{Random fields and probability spaces}
\label{Section:ProbaSpaces}

We introduce a  probability space $(\Omega,\mathcal{F},\mathbb{P})$, 
with $\Omega$ being the set of all elementary events, or space of outcomes,
and $\mathcal{F}$ a $\sigma$-algebra of all possible events, equipped with 
a probability measure $\mathbb{P}$. 
Random entropy solutions are random functions taking values in a function
space; to this end, let $(E,\mathcal{G},\mathbb{G})$ denote any measurable space. 
Then an $E$-valued random variable is any 
mapping $Y : \Omega \to E$ such that $\forall A \in \mathcal{G}$ the 
preimage $Y^{-1}(A) = \{\omega \in \Omega : Y(\omega) \in A\} \in \mathcal{F}$, i.e. such that $Y$ is a $\mathcal{G}$-measurable 
mapping from $\Omega$ into $E$. 

We confine ourselves to the case that $E$ is a complete metric space; 
then $(E,\mathcal{B}(E))$ equipped with a Borel $\sigma$-algebra $\mathcal{B}(E)$ 
is a measurable space. 
By definition, $E$-valued random variables $Y : \Omega \to E$ 
are $\big(E,\mathcal{B}(E)\big)$ measurable. 
Furthermore, if $E$ is a separable Banach space with norm $\|\circ\|_E$ 
and with topological dual $E^*$, then $\mathcal{B}(E)$ is the smallest 
$\sigma$-algebra of subsets of $E$ containing all sets
\[ 
\{x \in E : \varphi(x) < \alpha\}, \varphi \in E^*, \alpha \in \mathbb{R}\;.
\]
Hence, if $E$ is a separable Banach space, $Y : \Omega \to E$ is an 
$E$-valued random variable if and only if for every 
$\varphi \in E^*$, $\omega \mapsto \varphi\big(Y(\omega)\big) \in \mathbb{R}$ 
is an $\mathbb{R}$-valued random variable.
Moreover, there hold the 
following results on existence and uniqueness 
\cite{MishraSchwab2010}.

For a simple $E$-valued random variable $Y$ and for any $B \in \mathcal{F}$ we set
\begin{equation}
\label{Eq:pre1}
\int_B Y(\omega)\,\mathbb{P}(d\omega) 
= 
\int_B Y\,d\mathbb{P} 
= 
\sum\limits_{i=1}^N x_i\mathbb{P}(A_i\cap B).
\end{equation}

For such $Y(\omega)$ and all $B \in \mathcal{F}$ holds
\begin{equation}
\label{Eq:pre2}
\Big\|\int_B Y(\omega)\,\mathbb{P}(d\omega)\Big\|_E \leq \int_B \|Y(\omega)\|_E\,\mathbb{P}(d\omega). 
\end{equation}
For any random variable $Y : \Omega \to E$ which is Bochner integrable, 
there exists a sequence $\{Y_m\}_{m\in \IN}$ of simple random variables such that, 
for all $\omega \in \Omega, \|Y(\omega) - Y_m(\omega)\|_E \to 0$ as $m \to \infty$. 
Therefore \eqref{Eq:pre1} and \eqref{Eq:pre2} can be extended to any $E$-valued random variable. 
We denote the expectation of $Y$ by
\[\IE[Y] = \int_{\Omega} Y(\omega)\,\mathbb{P}(d\omega) = \lim\limits_{m\to\infty}\int_{\Omega}Y_m(\omega)\mathbb{P}(d\omega) \in E,
\]
and the variance of $Y$ is defined by
\[
\IV[Y] = \IE\big[(Y-\IE[Y])^2\big].
\]

Denote by $L^p(\Omega,\mathcal{F},\mathbb{P};E)$ 
for $1 \leq p \leq \infty$ the Bochner space of all 
$p$-summable, $E$-valued random variables $Y$ and 
equip it with the norm
\[
\|Y\|_{L^p(\Omega;E)} = \big(\IE[\|Y\|^p_E]\big)^{1/p} 
= 
\left(\int_{\Omega}\|Y(\omega)\|^p_E\,\mathbb{P}(d\omega)\right)^{1/p}.
\]

For $p=\infty$ we can denote by $L^{\infty}(\Omega,\mathcal{F},\mathbb{P};E)$ 
the set of all $E$-valued random variables which are essentially bounded and 
equip this space with the norm
\[
\|Y\|_{L^{\infty}(\Omega;E)} = \mathrm{ess}\sup\limits_{\omega\in\Omega}\|Y(\omega)\|_E.
\]

Consider now the balance law \eqref{eq:hyp_cons} and assume that the parameter $\mubold$ represents vector of real-valued real variables. Different uncertainty quantification (UQ) techniques can be applied to model the effects of this randomness in $\mubold$ on the solution $\ubold$.

\subsection{Monte Carlo method}

In this chapter, we restrict ourselves to the applications of ROM techniques to UQ problems in conjunction with the well-known Monte Carlo sampling method. We note, however, that the outlined ideas could be easily extended to more recent sampling methods such as Multi-Level Monte Carlo (MLMC) method, as well as Stochastic Collocation methods.

The idea of the Monte Carlo method consists in generating $M$ independent, identically distributed samples $\bar{\mubold}^i$ of the random variable $\mubold$, for $i = 1,\dots,M$, and calculating the corresponding deterministic approximate solutions $\bar{\ubold}^i$ of \eqref{eq:hyp_cons}. Then, the Monte Carlo estimate of the expected solution value $\IE[\ubold]$ at time $t$ and at point $x$ is given by
\begin{equation}
\label{eq:exp_mc}
E_M[\ubold(x,t)] = \dfrac1M \sum\limits_{i=1}^M \bar{\ubold}^i (x,t),
\end{equation}
and the variance can be computed according to the unbiased estimate
\begin{equation}
\label{eq:var_mc}
V_M[\ubold(x,t)] = \dfrac{1}{M-1}\sum\limits_{i=1}^M \big( \bar{\ubold}^i (x,t) - E_M \big)^2.
\end{equation}

\section{Numerical results}\label{ch:results}
In this chapter we will present our numerical results that illustrate the behavior of the RB methods in the case of nonlinear unsteady hyperbolic conservation laws in 1D and 2D with applications in UQ.
\subsection{Stochastic unsteady Burgers' equation in 1D with random data}\label{subsec:Burgers_eq} We consider here Burgers' equations with randomness in both flux and initial data
\begin{align}\label{eq:Burgers_1D}
	&\frac{\partial u}{\partial t}+\frac{\partial f(u,w)}{\partial x} =0, ~x\in[0,\pi], ~w\in \Omega,  \\
	&u_0(x,w) =u_0(x,Y_1(w),Y_2(w)),
\end{align}
defined on $D=[0,\pi]\subset \R, ~t>0$ with periodic boundary conditions, the nonlinear flux is given as:
\begin{equation}\label{eq:Burgers_flux}
f(u,w)=f(u,Y_3(w))=Y_3(w)f(u)=Y_3(w)\frac{u^2}{2}
\end{equation}
and the initial condition is given by:
\begin{equation}\label{eq:IC_Burgers_eq}
u_0(x,Y_1(w),Y_2(w))=|\sin(2x+Y_1(w))|+0.1Y_2(w),
\end{equation}
where $y_j=Y_j(w),~ j=1,2,3, ~w \in \Omega  $ and $Y_j$ is a random variable which takes values in the domain $\PP\subset \R^q$ of the parametrized probability space.

The PDE is discretized by an upwind first order finite volume scheme. We used an uniform mesh $\{ x_{i-1/2} \}_{i=1}^{N_h+1}$, resulting in a HDM of dimension $N_h=10^3$, with the CFL condition of 0.318, $K=159$ time iterations, final time $t^K=0.159$ and time step of 0.001. In this first example, we will use a finite volume approach, in the RD context, since it can be rewritten in this formulation thanks to \cite{abgrall2017}. With $x_{i-1/2}$ defining the points of the grid, we define the cells $T_{i}=[ x_{i-1/2}, x_{i+1/2} ]$ and we consider constant approximation over each cell $u_{i}$. The scheme will then read $u^{k+1}_{i}=u^{k}_{i} -\frac{\Delta t}{\Delta x} \left( f_{i+1/2} - f_{i-1/2}\right)$  . We are using the numerical Roe fluxes $f$ defined at the cell interface as:
\begin{equation}\label{eq:Roe_flux_general}
f_{i+1/2}=f(u_L,u_R)=\frac{1}{2}\Big[f(u_L)+f(u_R)-|a(u_L,u_R)|(u_R-u_L) \Big],
\end{equation}
where $u_L=u_i$ and $u_R=u_{i+1}$. The Rankine-Hugoniot velocity is
$$a(u_L,u_R)=\frac{f(u_L)-f(u_R)}{u_L-u_R}.$$
This numerical flux choice has the purpose of linearizing the flux $f$ around the cell interface and then using an upwind flux, which has the role of an entropy fix. For Burgers' equations, the Roe flux including the randomness $Y_3(w)$ writes
\begin{equation}\label{eq:Roe_flux_Burgers}
f(u_L,u_R)=\frac{1}{4}Y_3(w)\Big[u_L^2+u_R^2-|u_L^2+u_R^2|(u_R^2-u_L^2) \Big].
\end{equation}
We consider now two cases: the first one which consists only in one randomness in the initial data and the second case which contains randomness in the flux and in the initial condition.

\subsubsection{Stochastic unsteady Burgers' equation with random initial data}  
\begin{figure}[h!]
	\begin{center}
		{\includegraphics[width=0.55\textwidth,clip=]{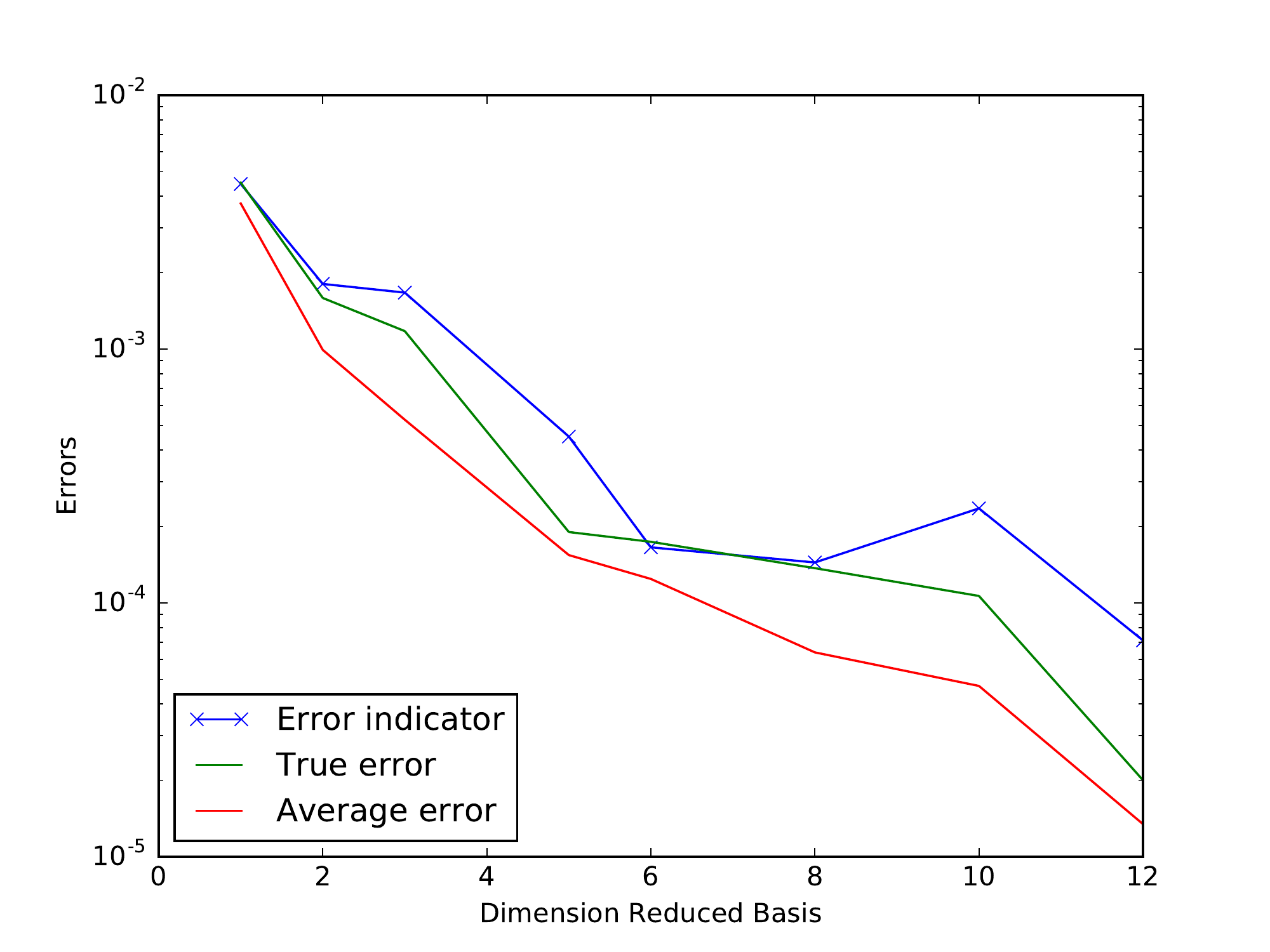}}
	\end{center}
	\caption{\label{fig:error_Burgers_1_param} The error decrease during basis extension with growing RB size for Burgers' equation with one random data}
\end{figure}

In this case, we consider as deterministic $Y_2(w)=Y_3(w)=1,\, \forall w \in \Omega$, while $Y_1(w) \sim \UU [0.4,0.5] $ is the only random variable. In the greedy procedure we sampled the training set using an uniform grid on the parameter domain $D_y=[0.4, 0.5]$. We have not used the PODEIM--Greedy algorithm in this test case (the EIM is performed before the POD--Greedy), because the error of the greedy procedure was naturally decreasing without oscillations. The tolerance set for the EIM procedure was $10^{-6}$ and for the greedy algorithm was $10^{-4}$. What we get from offline phase is an EIM space with 61 functions and a RB space of dimension 12 (see Figure \ref{fig:error_Burgers_1_param}). 

\begin{figure}[h!]
	\begin{center}
		{\includegraphics[width=0.5\textwidth,clip=]{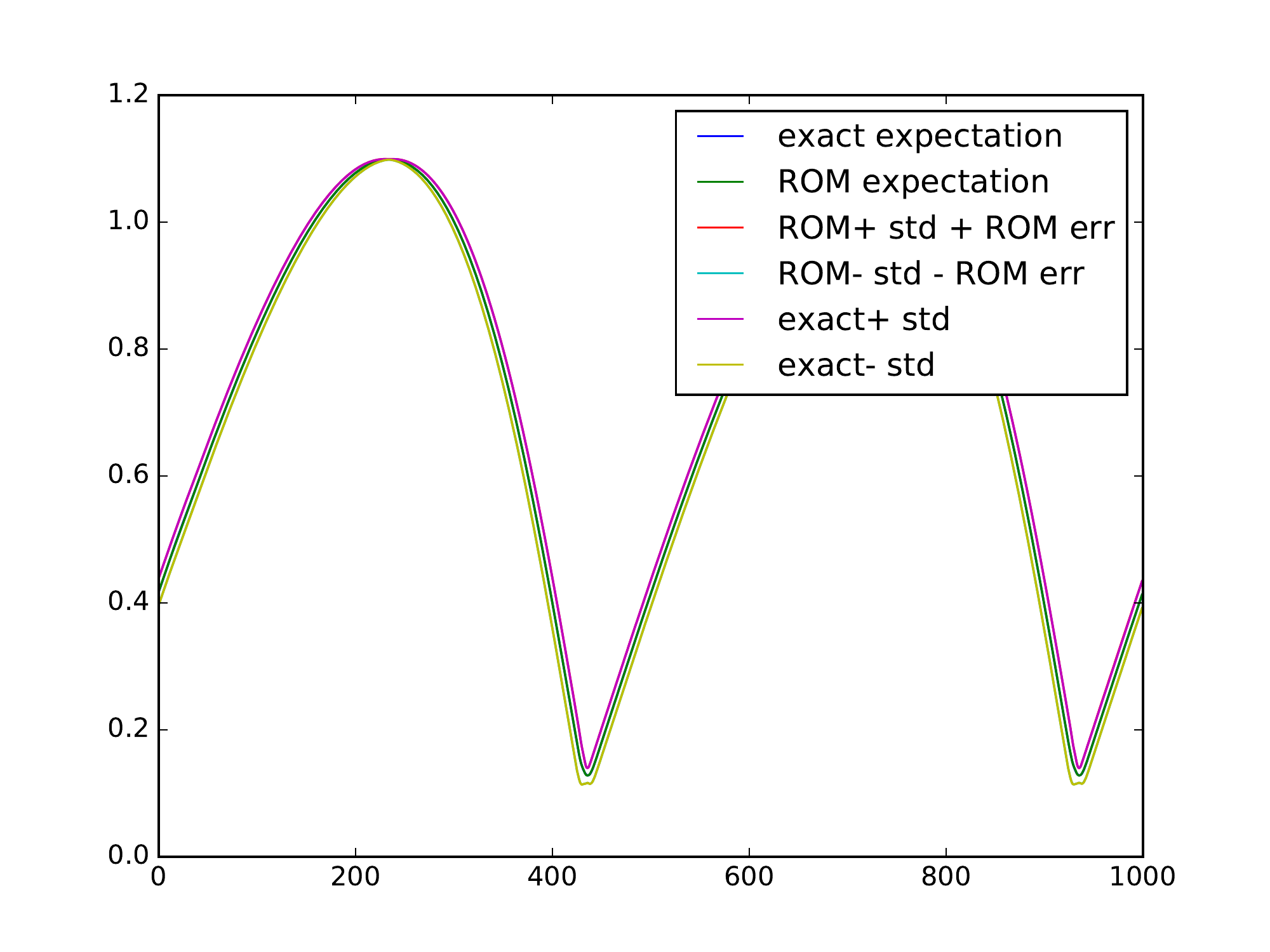}}
	\end{center}
	\caption{\label{fig:mean_Burgers_1_param} Solution mean and the mean plus/minus the standard deviation for both the reduced and the high-fidelity problem in the case of Burgers' equation with one random data }
\end{figure}

\begin{figure}[h!]
	\begin{center}
		{\includegraphics[width=0.5\textwidth,clip=]{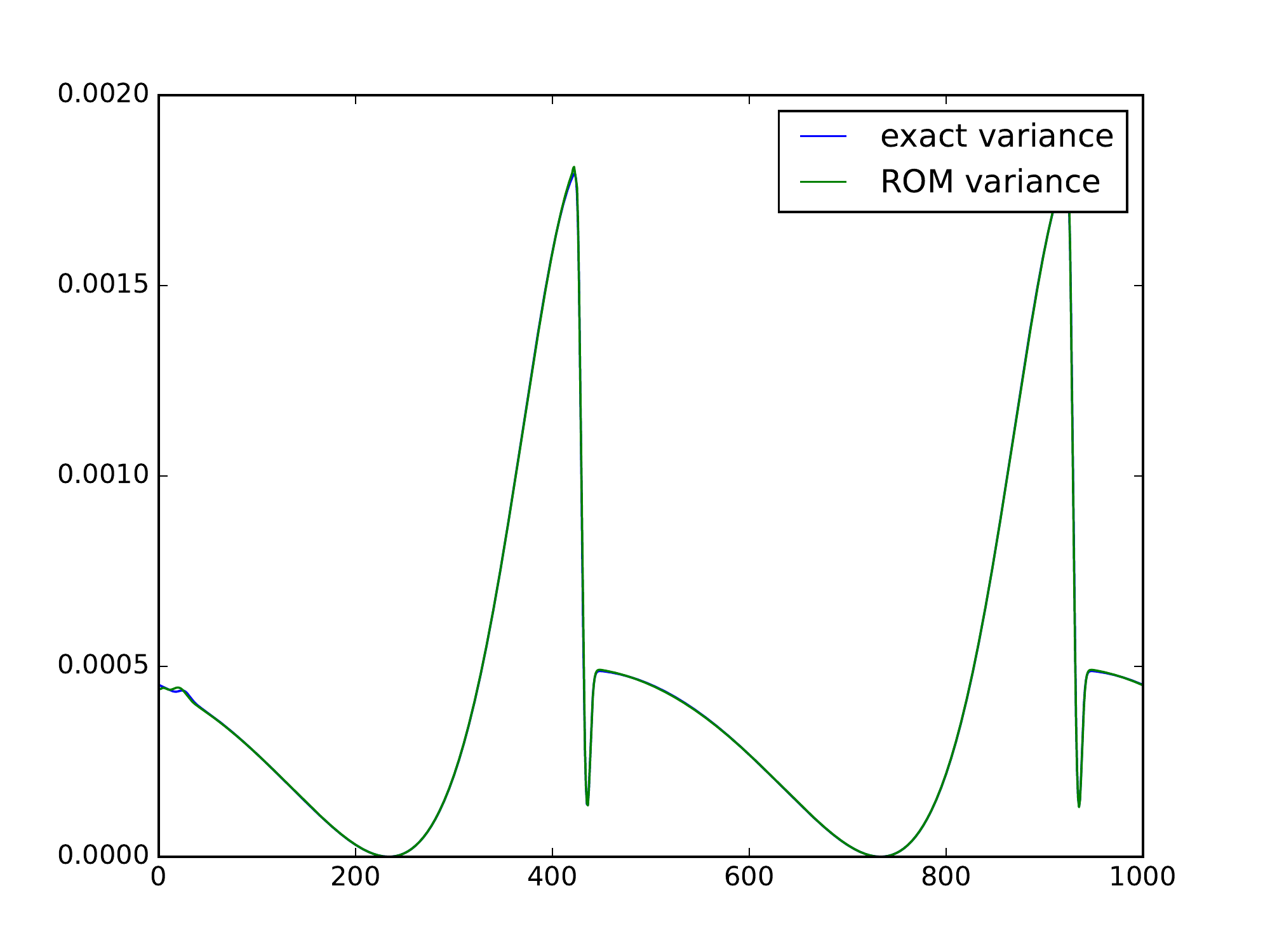}}
	\end{center}
	\caption{\label{fig:var_Burgers_1_param} Variance for the reduced and the high-fidelity problem in the case of Burgers' equation with one random data }
\end{figure}

For the online phase, we want to compute some statistical moments with arbitrary probability distributions of the uncertainty, such as the solution mean and the variance, as well as the solution mean plus/minus the standard deviation of the random variable $u_h^{K}(w)$. This UQ analysis is performed using a set with 100 elements in the parameter domain $D_y=[0.4,0.5]$, which were generated by a random Monte Carlo method. The advantage of performing an UQ analysis after a RB procedure is that the computational time for a single reduced solution will be much lower than the high fidelity one, the solution accuracy being comparable (see Figure \ref{fig:mean_Burgers_1_param}, \ref{fig:var_Burgers_1_param}). Indeed, the average computational time for one high fidelity solution is of 1.2551 seconds, while the reduced solution takes only 0.17118 seconds, the percentage of the saved time being then of 86\%. \footnote{The computations are performed with a Intel(R) Xeon(R) CPU E7-2850  @ 2.00GHz}



\subsubsection{Stochastic unsteady Burgers' equation with random flux and initial data}  
\begin{figure}[h!]
	\begin{center}
		{\includegraphics[width=0.55\textwidth,clip=]{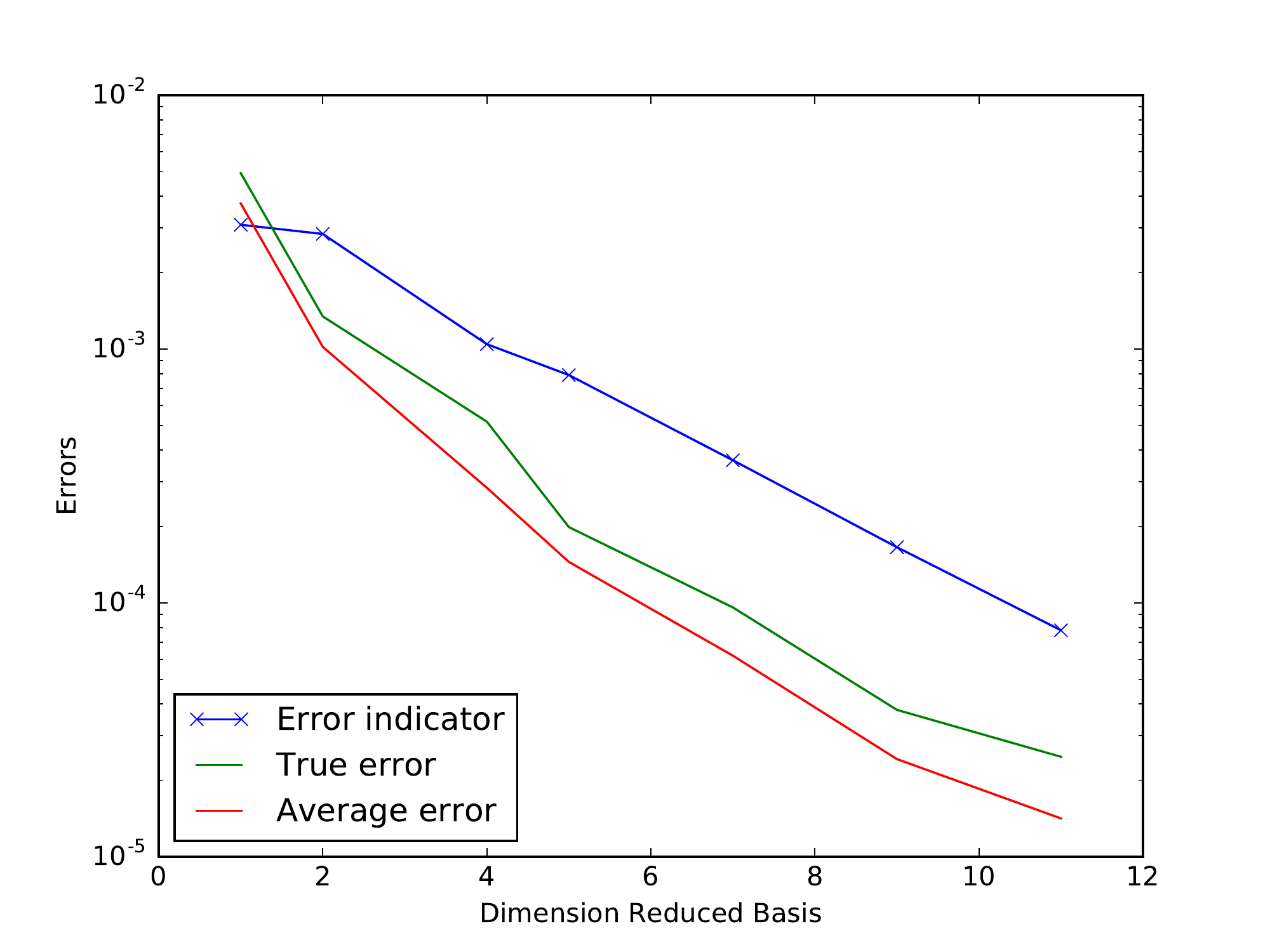}}
	\end{center}
	\caption{\label{fig:error_Burgers_3_param} The error decrease during basis extension with growing RB size for Burgers' equation with random flux and random initial condition}
\end{figure}

Consider now the case of Burgers' equation with randomness in both flux and initial condition, namely $Y_3(w)$, respectively $Y_1(w)$ and $Y_2(w)$. Let us define $Y_1\sim \UU [0.4,0.5],\,Y_2\sim \UU [1,1.2] ,\,Y_3\sim \UU [0.9,1.1]$. In the greedy procedure we sampled the training set using an uniform three-dimensional grid on the parameter domain $D_y=[0.4, 0.5] \times [1,1.2] \times [0.9,1.1]$. We are using the same tolerances for the construction of the EIM space and of the RB as in the previous test case and without using any PODEI algorithm, we obtain an EIM space with 48 functions and an RB space of dimension 11 (see Figure \ref{fig:error_Burgers_3_param}).

\begin{figure}[h!]
	\begin{center}
		{\includegraphics[width=0.5\textwidth,clip=]{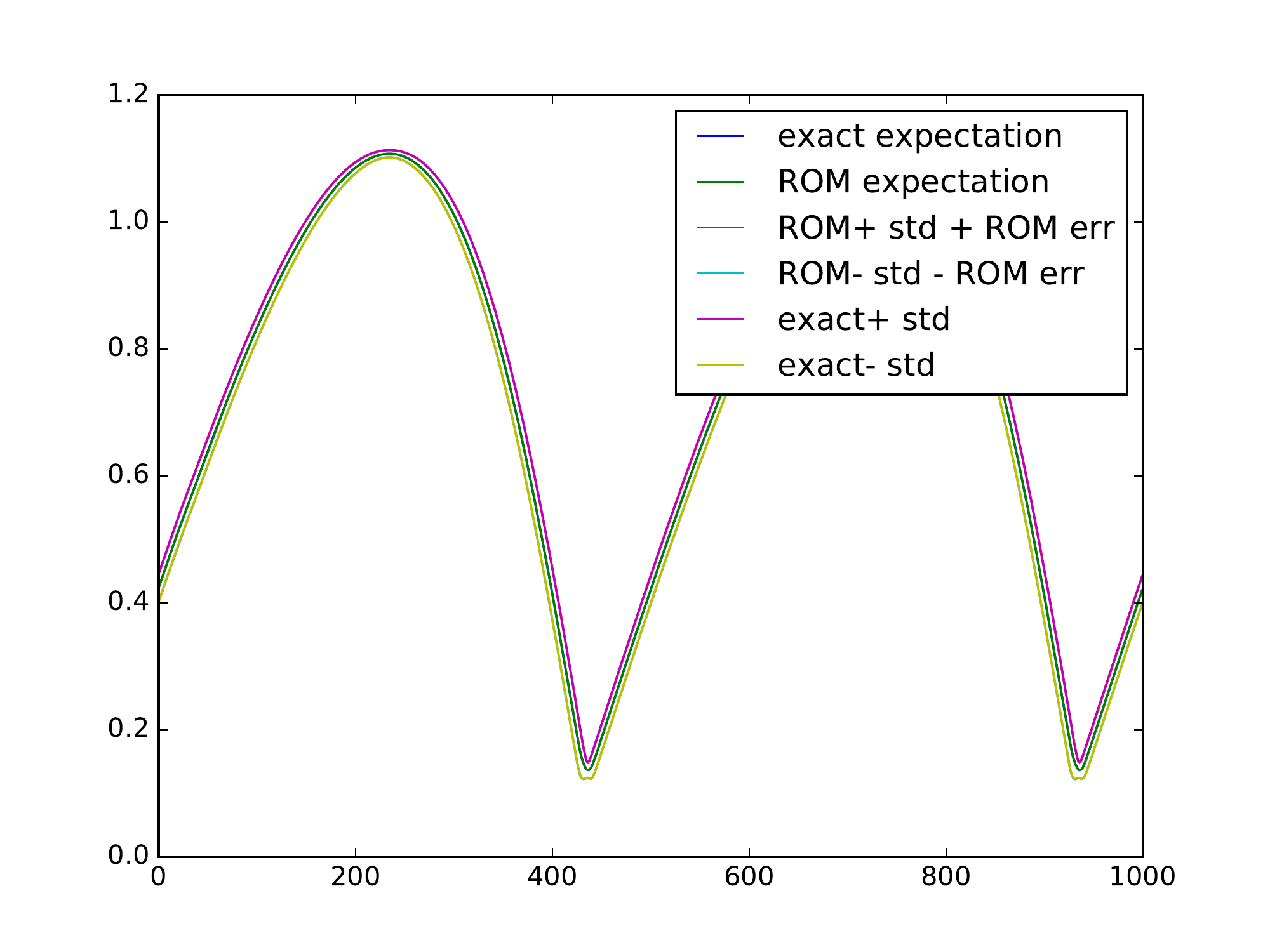}}
	\end{center}
	\caption{\label{fig:mean_Burgers_3_param} Solution mean and the mean plus/minus the standard deviation for both the reduced and the high-fidelity problem in the case of Burgers' equation with random flux and random initial condition}
\end{figure}

\begin{figure}[h!]
	\begin{center}
		{\includegraphics[width=0.5\textwidth,clip=]{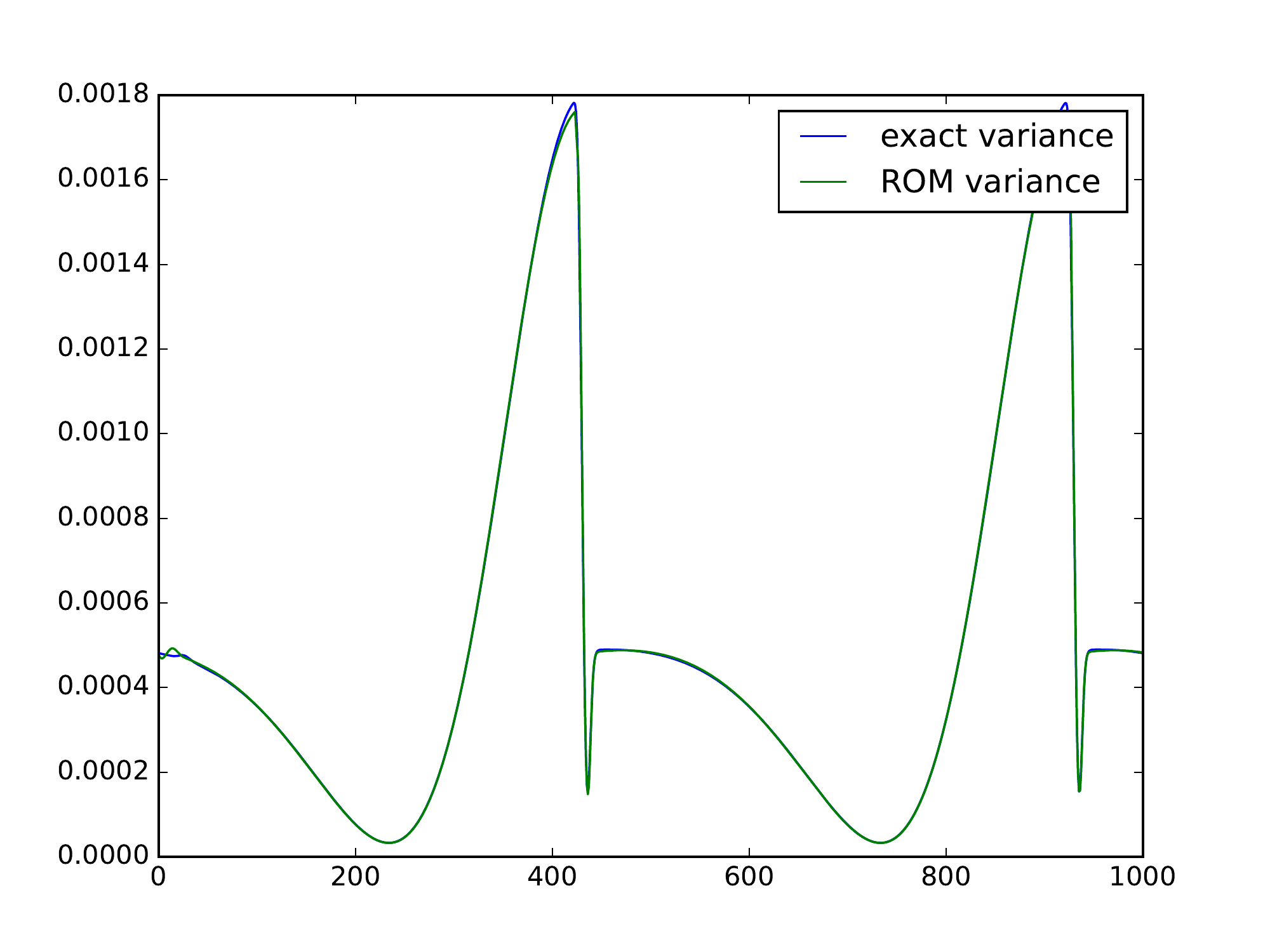}}
	\end{center}
	\caption{\label{fig:var_Burgers_3_param} Variance for the reduced and the high-fidelity problem in the case of Burgers' equation with random flux and random initial condition}
\end{figure}

In the online phase, the UQ analysis is performed using a set with 125 elements in the parameter domain $D_y=[0.4, 0.5] \times [1,1.2] \times [0.9,1.1]$, which were generated by a random Monte Carlo method.  Comparing again the solution mean and the variance, as well as the solution mean plus/minus the standard deviation of a random variable $u_h^{K}(w)$ in the case of the reduced problem and the high fidelity one (see Figure \ref{fig:mean_Burgers_3_param}, \ref{fig:var_Burgers_3_param}), we obtain a computational saving time of 88\%. Indeed, the average computational time for one high fidelity solution is of 1.2143 seconds, while the reduced solution takes only 0.14472 seconds.



\subsection{Stochastic Euler equations in 1D with random data}\label{subsec:Euler_1D}
We consider the parametrized Euler equations
\begin{align}\label{eq:Euler_1D}
	\frac{\partial \ubold}{\partial t}+\frac{\partial \fbold(\ubold,w)}{\partial x}=0, ~x \in [-1,1] \\
	\ubold_0(x,w)=\ubold_0(x,Y_1(w))
\end{align} 
with $y_j=Y_j(w), ~ j=1,2 ~w \in \Omega$ and 
$$\ubold=(\rho,\rho u, E )^{T}, ~ \fbold=(\rho,\rho u^2+p,\rho u(E+p))^T, ~p=(\gamma-1)(E-\frac{1}{2}\rho u^2).$$
We also assume the randomness in the adiabatic constant, $\gamma=Y_2(w)$, and therefore the flux is parameter dependent: $$\fbold(\ubold,w)=\fbold(\ubold,Y_2(w)).$$ 

We consider again two cases: the first one when we have randomness only in the initial data and the second case when we have randomness in the initial data and also in the specific heat ratio $\gamma$.

\subsubsection{Stochastic Euler equations in 1D with random initial data} 
For this smooth test case, we consider the following random initial condition:
$$\ubold_0(x,Y_1(w))=\Big(2+\sin(30Y_1(w))\sin(\pi(x-1)+Y_1(w)) ,0,(2+\sin(30Y_1(w))\sin(\pi(x-1)+Y_1(w)))^{\gamma}   \Big).$$
We set the value of the specific heat to $\gamma=Y_2(w)=1.4$ and we construct $Y_1(w)$ using a random Monte Carlo sampling method in the interval $D_y=[0.4,0.5]$, resulting in a set with 100 elements. The PDE is discretized by a first order finite volume scheme with MUSCL extrapolation on the characteristic variables and minmod limiter on all waves and the resulting HDM is of dimension $N_h=1200$ using $K=200$ time iterations of step $0.001$, final time $t^K=0.2$ and the space step of $0.001667$. 

\begin{figure}[h!]
	\begin{center}
		{\includegraphics[width=0.55\textwidth,clip=]{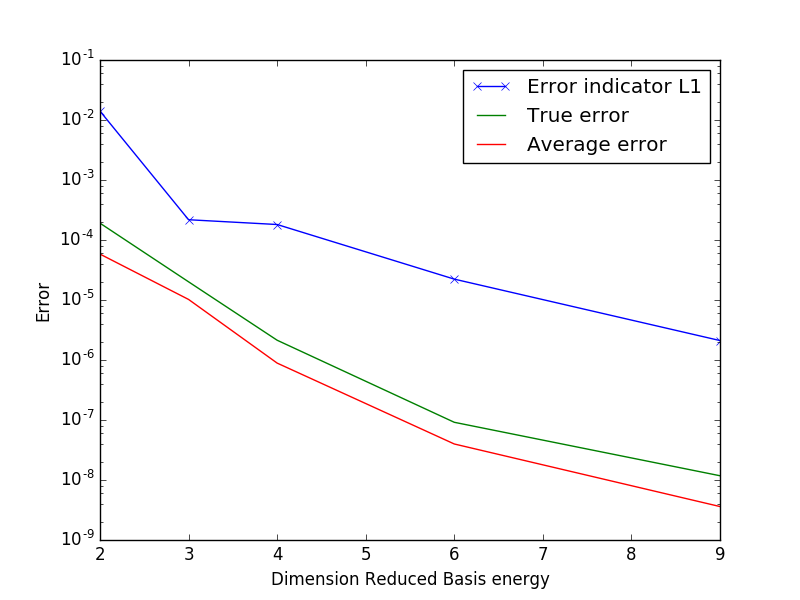}}
	\end{center}
	\caption{\label{fig:error_Euler_1_param} The error decrease during basis extension with growing RB size for the total energy component of Euler equation with one random data}
\end{figure}

In the offline step, the tolerance set for the greedy algorithm is $5\cdot 10^{-6}$ and we are using a PODEIM--Greedy algorithm generating an EIM space with $(10,11,10)$ basis and a RB space of dimension $(9,10,9)$ in each component, namely in density, momentum and total energy (see Figure \ref{fig:error_Euler_1_param} for the total energy). The PODEIM--Greedy algorithm helps us to avoid the unstable behaviour of the scheme. Indeed, if the accuracy of the empirical interpolation is not enough with respect to the accuracy of the RB space, namely we see an increment in the error, then we discard the newly computed RB functions. This will lead to an automatic control of the correlation between the dimension of the EIM space $N_{EIM}$ and the one of the RB space $N$, as seen also for this test case.

\begin{figure}[h!]
	\centering
	\subfigure[]{\includegraphics[width=0.49\textwidth,clip=]{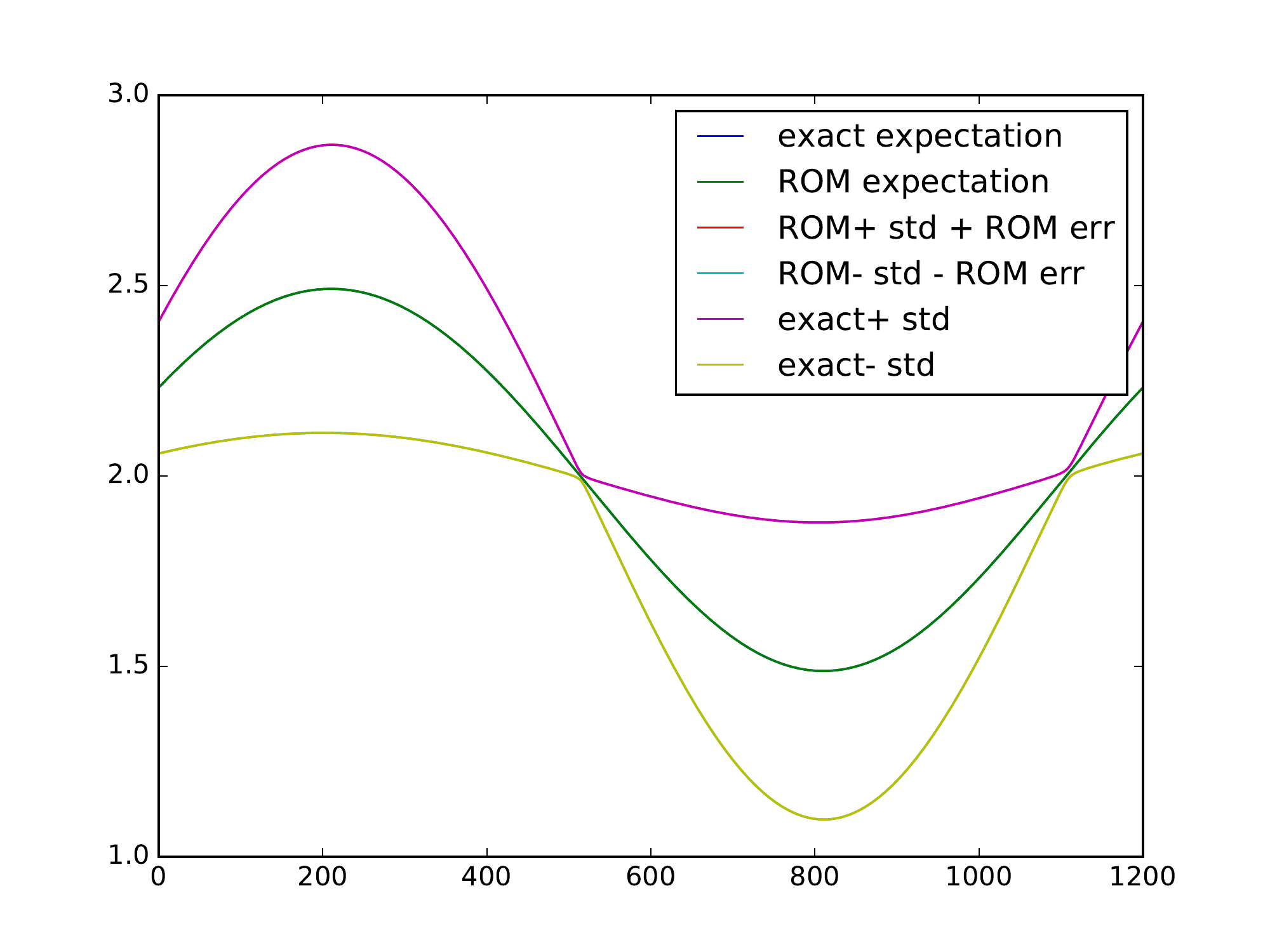}}
	\subfigure[]{\includegraphics[width=0.49\textwidth,clip=]{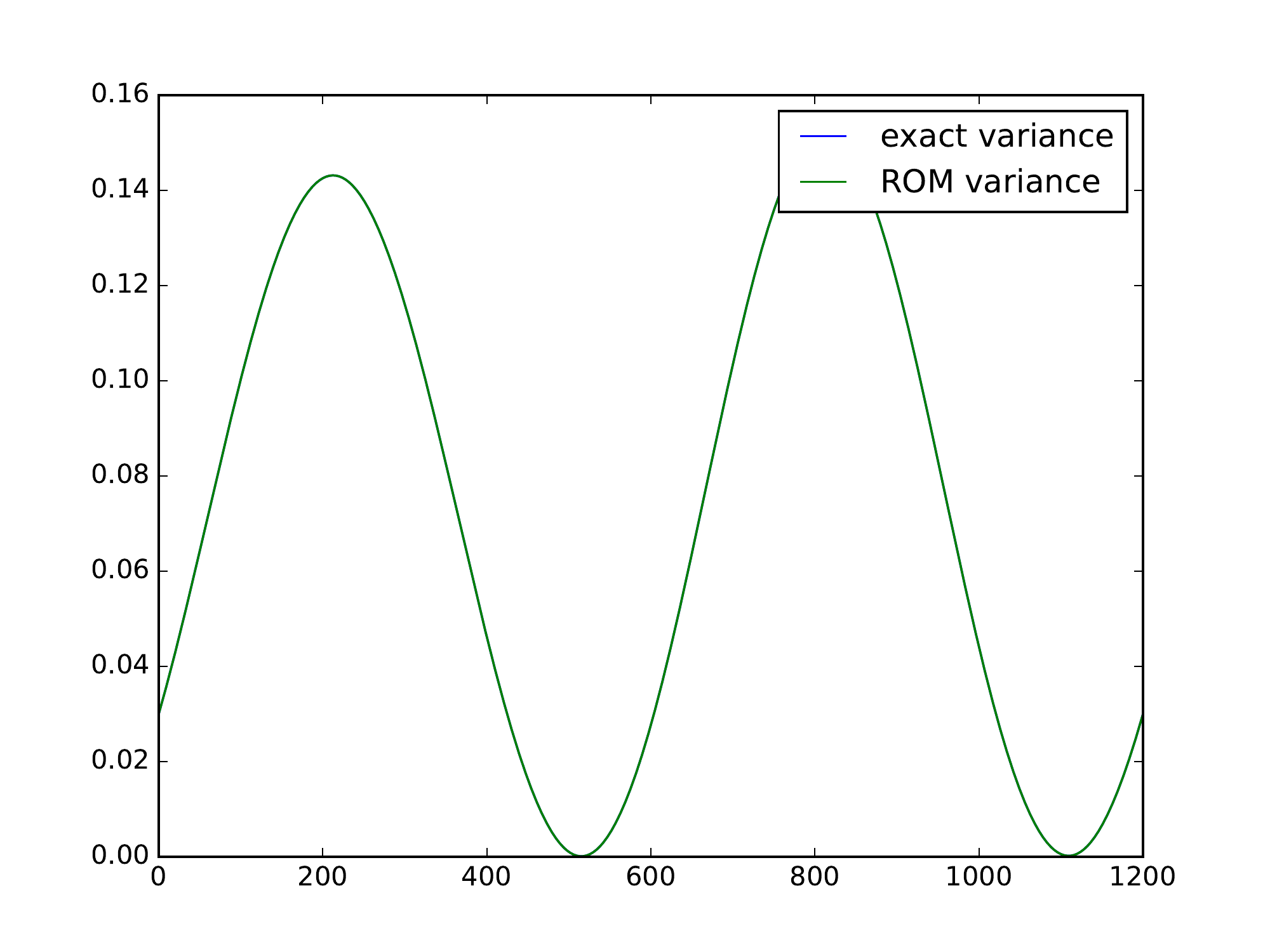}}
	\caption{\label{fig:mean_Euler_1d_1_param_rho} Solution mean, the mean plus/minus the standard deviation and the variance for both the reduced and the high-fidelity problem in the case of Euler equation with random initial condition for density}
\end{figure}

\begin{figure}[h!]
	\centering
	\subfigure[]{\includegraphics[width=0.49\textwidth,clip=]{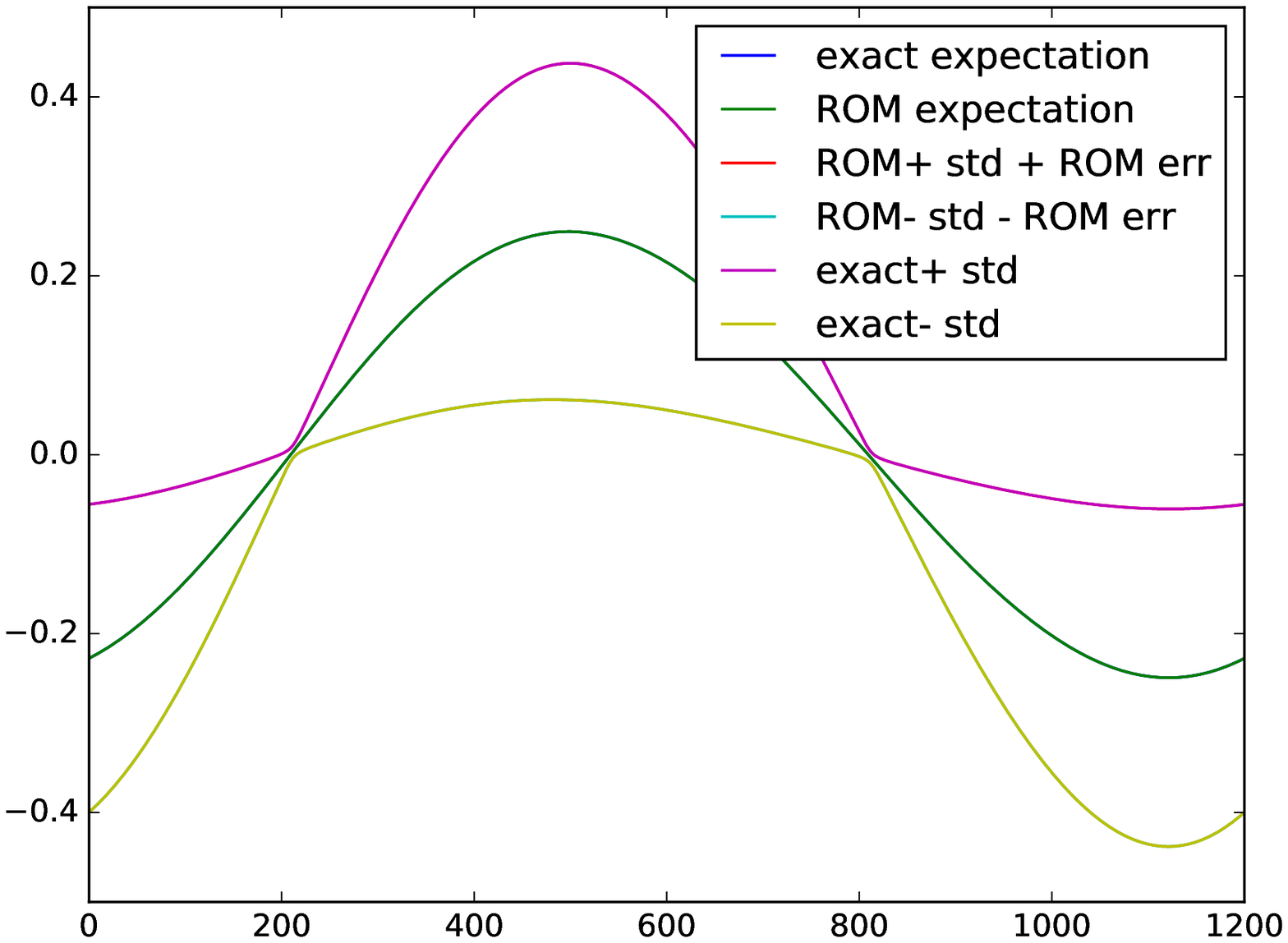}}
	\subfigure[]{\includegraphics[width=0.49\textwidth,clip=]{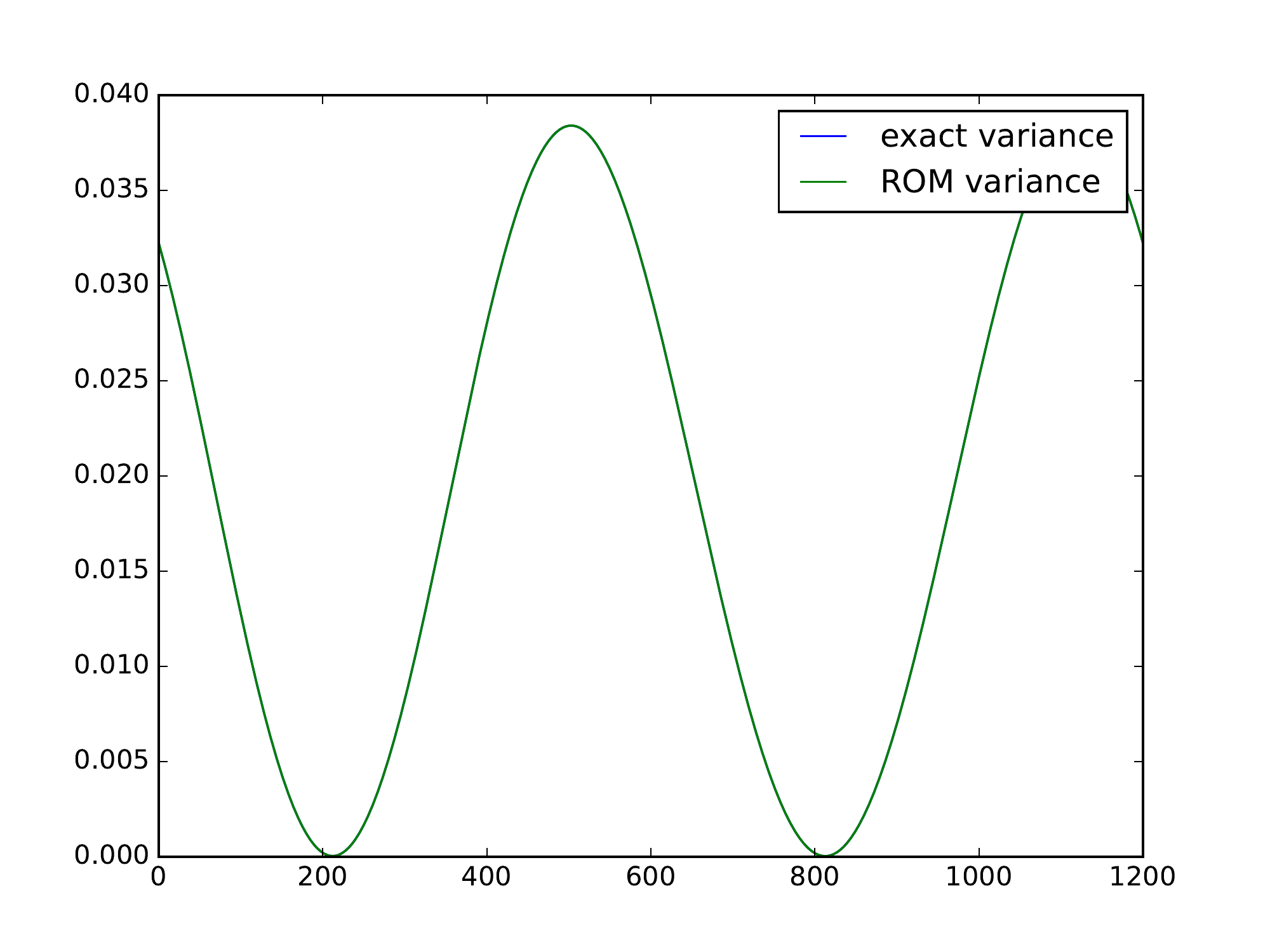}}
	\caption{\label{fig:mean_Euler_1d_1_param_mom} Solution mean, the mean plus/minus the standard deviation and the variance for both the reduced and the high-fidelity problem in the case of Euler equation with random initial condition for momentum}
\end{figure}

\begin{figure}[h!]
	\centering
	\subfigure[]{\includegraphics[width=0.49\textwidth,clip=]{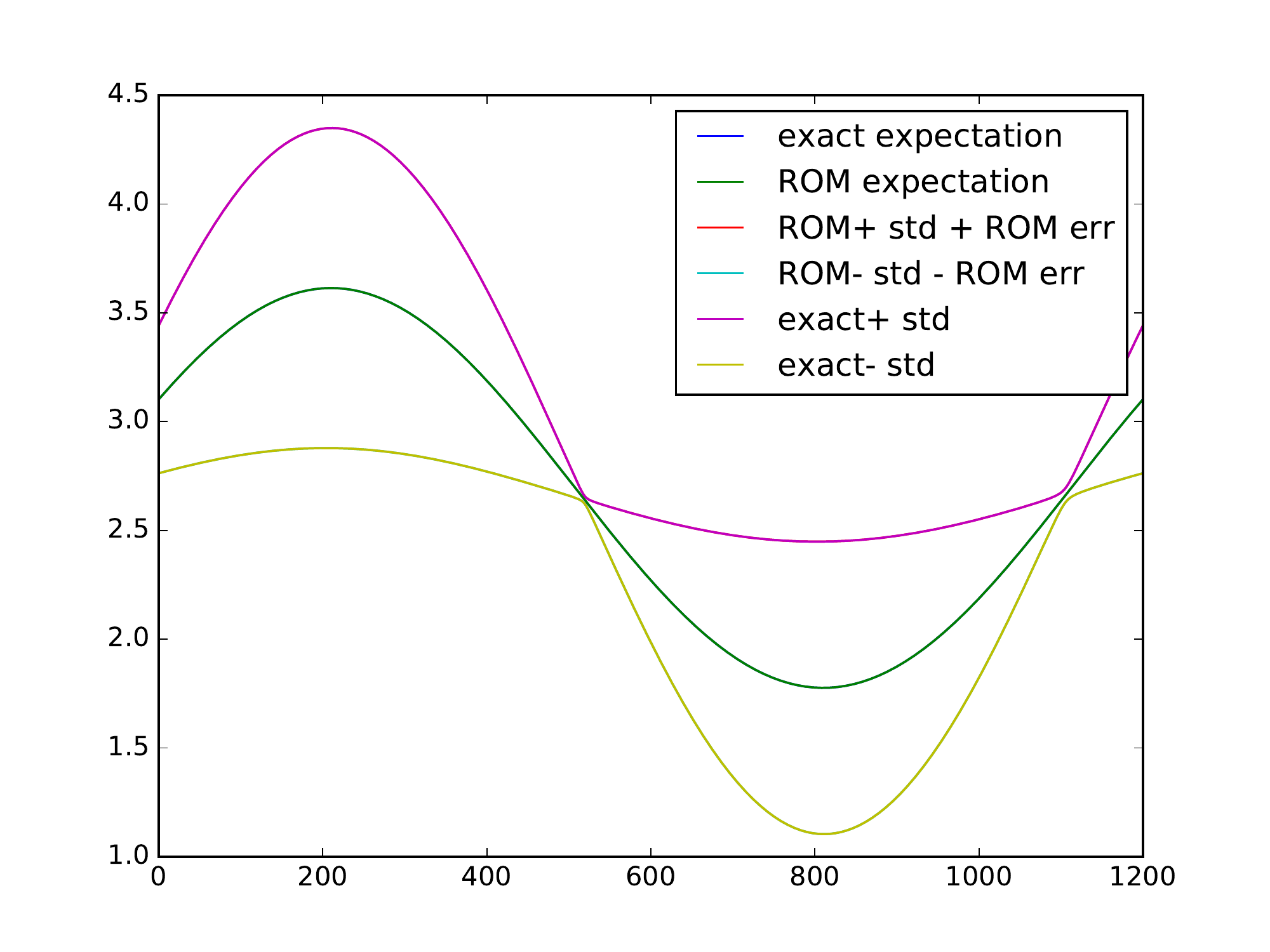}}
	\subfigure[]{\includegraphics[width=0.49\textwidth,clip=]{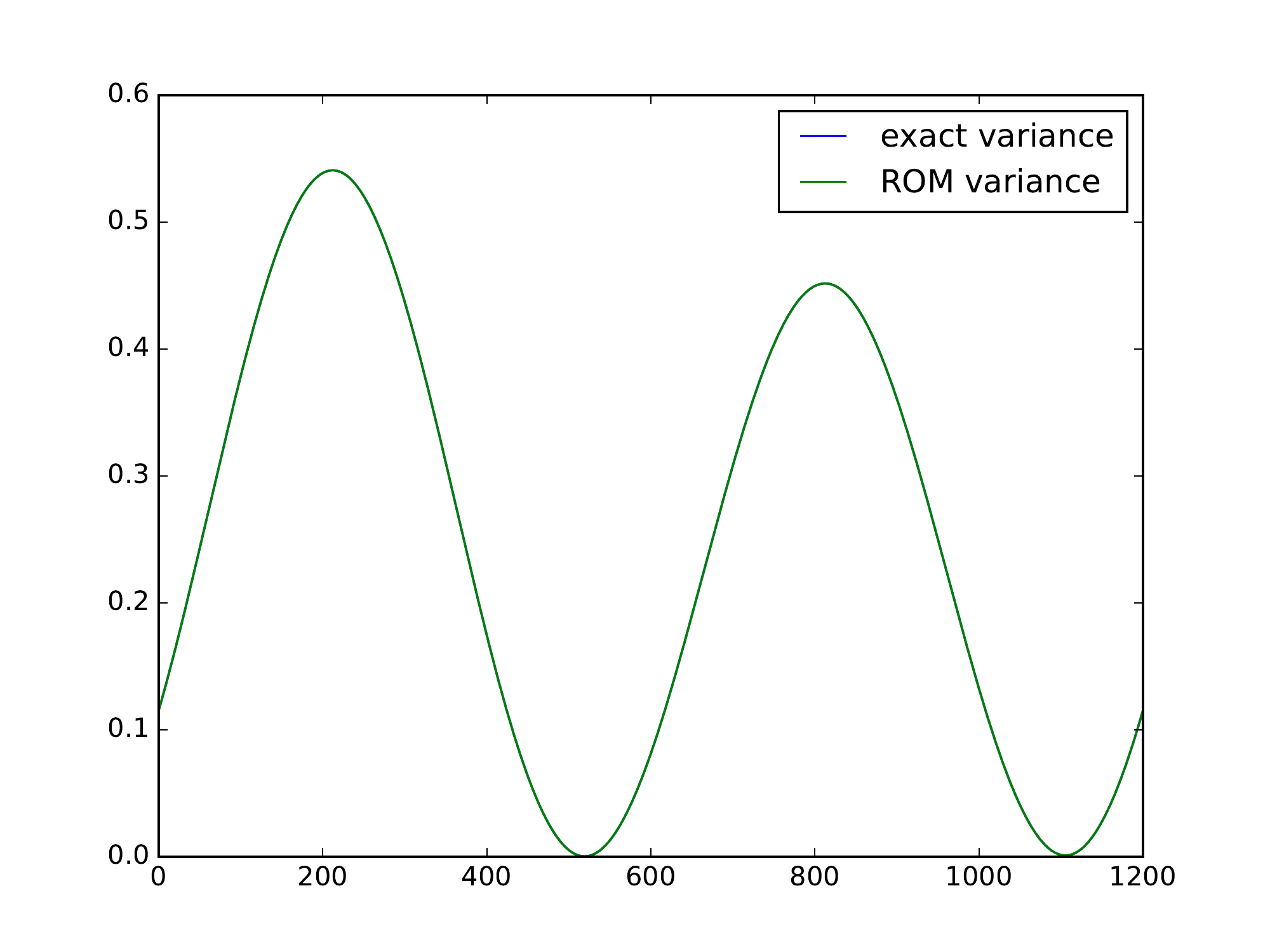}}
	\caption{\label{fig:mean_Euler_1d_1_param_E} Solution mean, the mean plus/minus the standard deviation and the variance for both the reduced and the high-fidelity problem in the case of Euler equation with random initial condition for the total energy}
\end{figure}

In the online phase, the UQ analysis is performed using a set with 100 samples in the parameter domain $D_y=[0.4, 0.5]$, which were generated by a random Monte Carlo method. Comparing again the solution mean and the variance, as well as the solution mean plus/minus the standard deviation of a random variable $\utruth{K}(w)$ in the case of the reduced problem and the high fidelity one (see Figures \ref{fig:mean_Euler_1d_1_param_rho}, \ref{fig:mean_Euler_1d_1_param_mom}, \ref{fig:mean_Euler_1d_1_param_E}), we obtain a computational saving time of 89\%. For a better visualization, we plot each component of the solution independently.  Indeed, the average computational time for one high fidelity solution is of 28.107 seconds, while the reduced solution takes only 3.2133 seconds.

\subsubsection{Stochastic Sod's shock tube problem in 1D with random initial data and random flux}
Consider now the Riemann problem for the one-dimensional Euler equations (\ref{eq:Euler_1D}) with the following initial data set in primitive variables:
\begin{equation*}
	\wbold_0(x,w)=(\rho_0(x,w),u_0(x,w),p_0(x,w))^T=\begin{cases}
		(1,0,1), & \mbox{if } x<0 \\ 
		(0.125+Y_1(w),0,0.1), & \mbox{if } x>0 .
	\end{cases}
\end{equation*}
In this test case, we have randomness in both flux and initial condition, namely the adiabatic constant $\gamma=Y_2(w)$, respectively $Y_1(w)$. We construct the random variables $Y_1(w),Y_2(w)$ using a random Monte Carlo sampling method in the interval $D_y=[-0.02,0.02] \times [1.4,1.5] $, resulting in a set with 100 samples. The PDE is discretized by a first order finite volume scheme with MUSCL extrapolation on the characteristic variables and minmod limiter on all waves and the resulting HDM is of dimension $N_h=1200$ using $K=320$ time iterations of step $0.0005$, final time $t^K=0.16$ and the space step of $0.001667$. 

\begin{figure}[h!]
	\begin{center}
		{\includegraphics[width=0.55\textwidth,clip=]{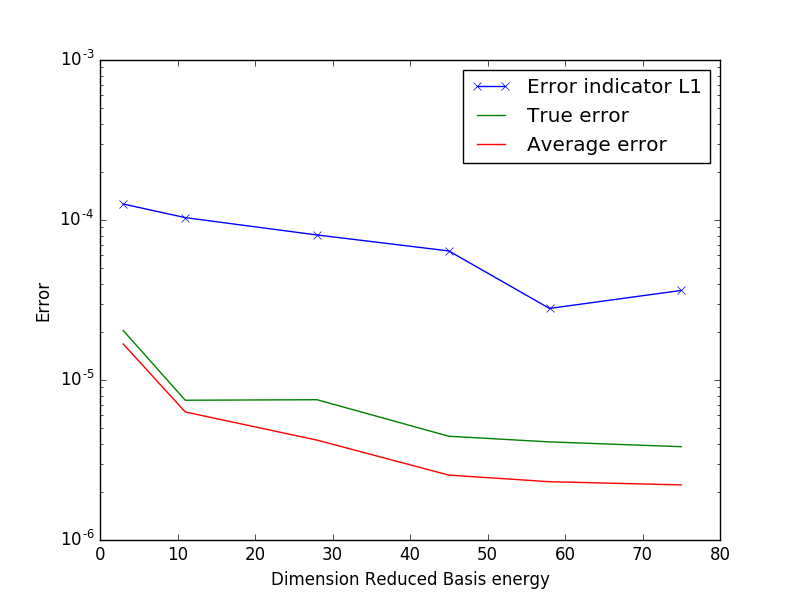}}
	\end{center}
	\caption{\label{fig:error_Euler_2_param} The error decrease during basis extension with growing RB size for the total energy component of Euler equation with one random data}
\end{figure}

In the offline step, the tolerance set for the greedy algorithm is $4\cdot10^{-6}$ and we are using a PODEI algorithm generating an EIM space with $ (68,83,89)$ basis and a RB space of dimension $(60,88,75)$ in each component, namely in density, momentum and total energy (see Figure \ref{fig:error_Euler_2_param} for the total energy).

\begin{figure}[h!]
	\centering
	\subfigure[]{\includegraphics[width=0.49\textwidth,clip=]{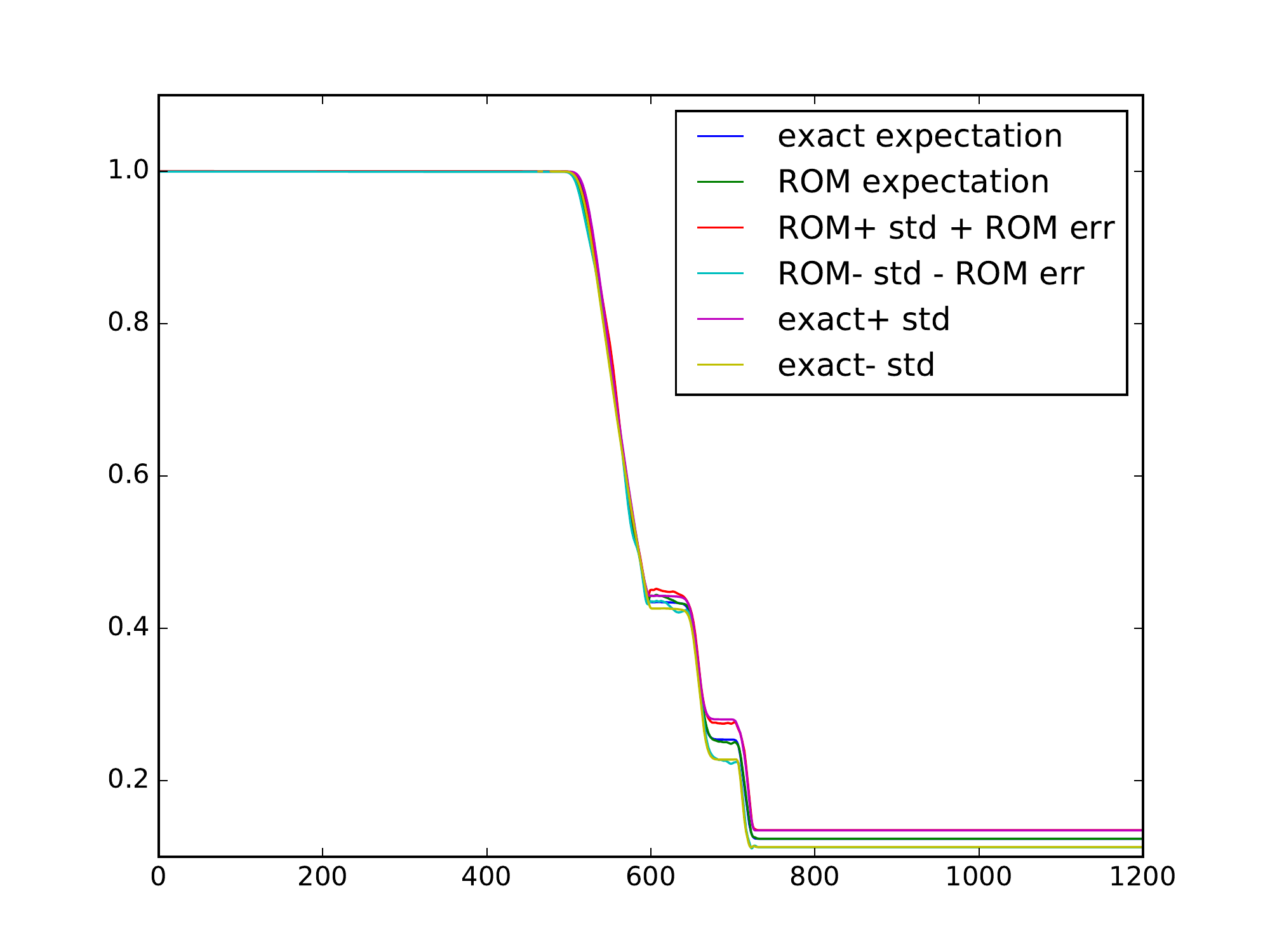}}
	\subfigure[]{\includegraphics[width=0.49\textwidth,clip=]{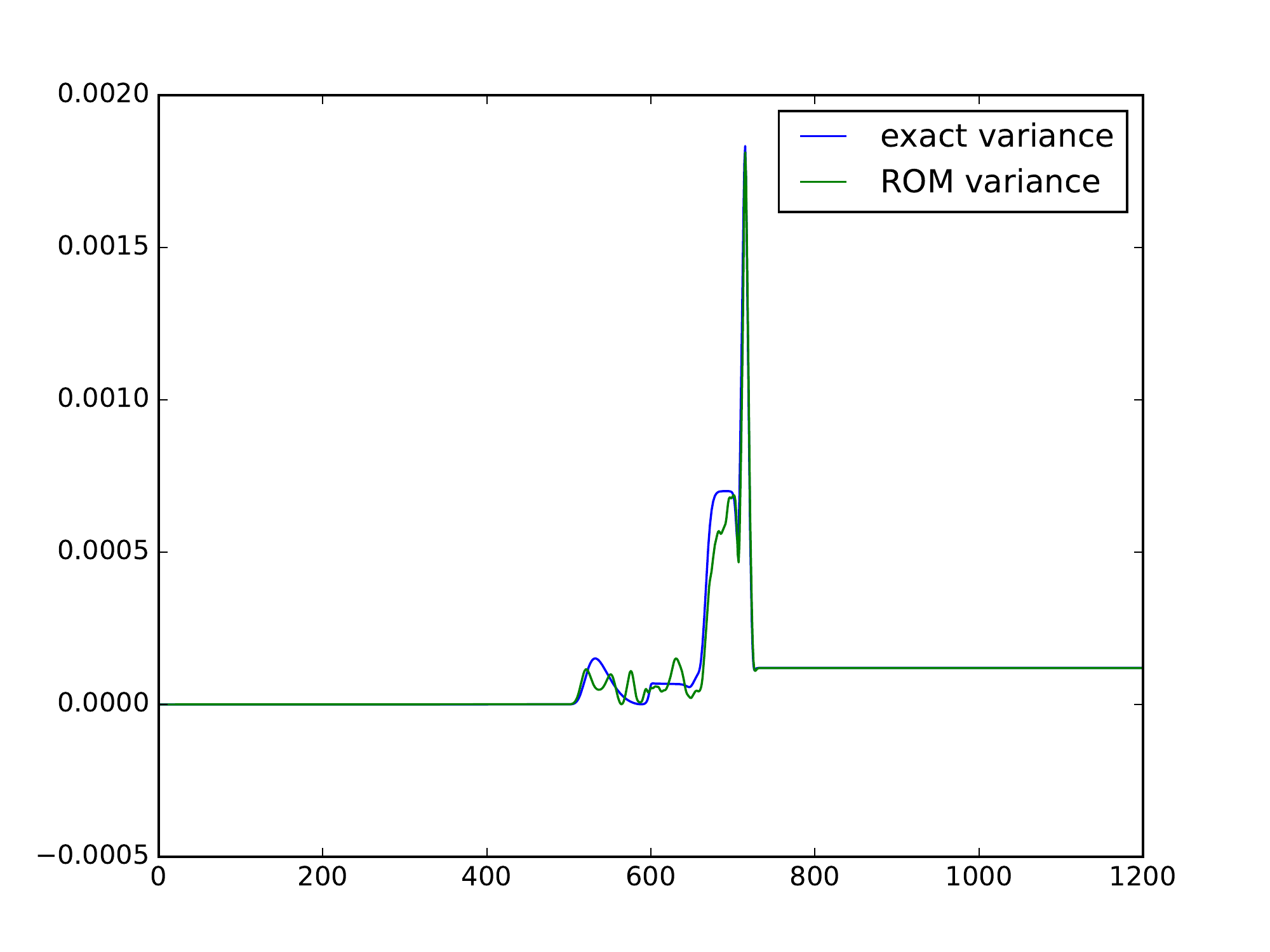}}
	\caption{\label{fig:mean_Euler_1d_2_param_rho} Solution mean, the mean plus/minus the standard deviation and the variance for both the reduced and the high-fidelity problem in the case of Euler equation with random initial condition and random flux for density}
\end{figure}

\begin{figure}[h!]
	\centering
	\subfigure[]{\includegraphics[width=0.49\textwidth,clip=]{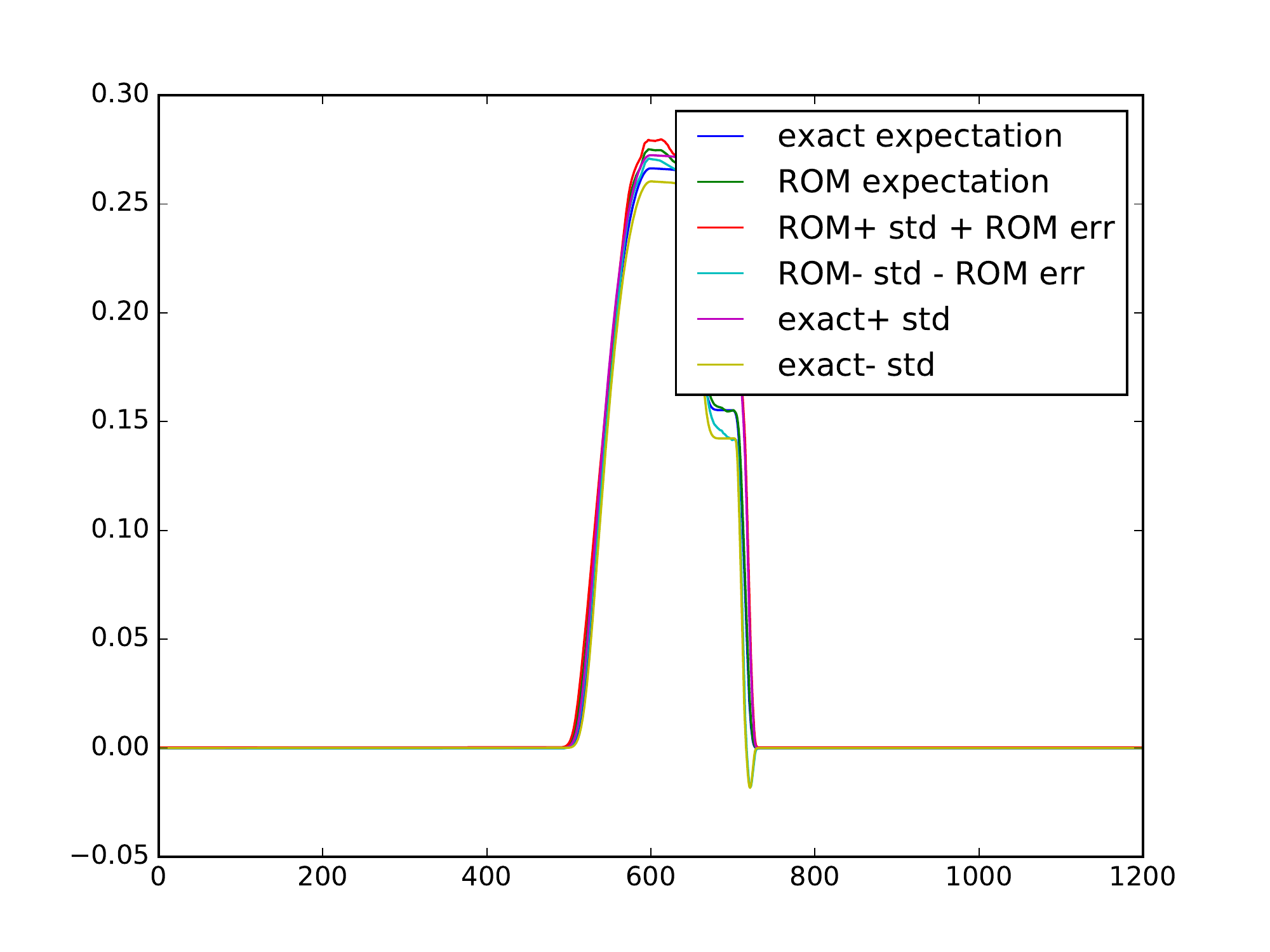}}
	\subfigure[]{\includegraphics[width=0.49\textwidth,clip=]{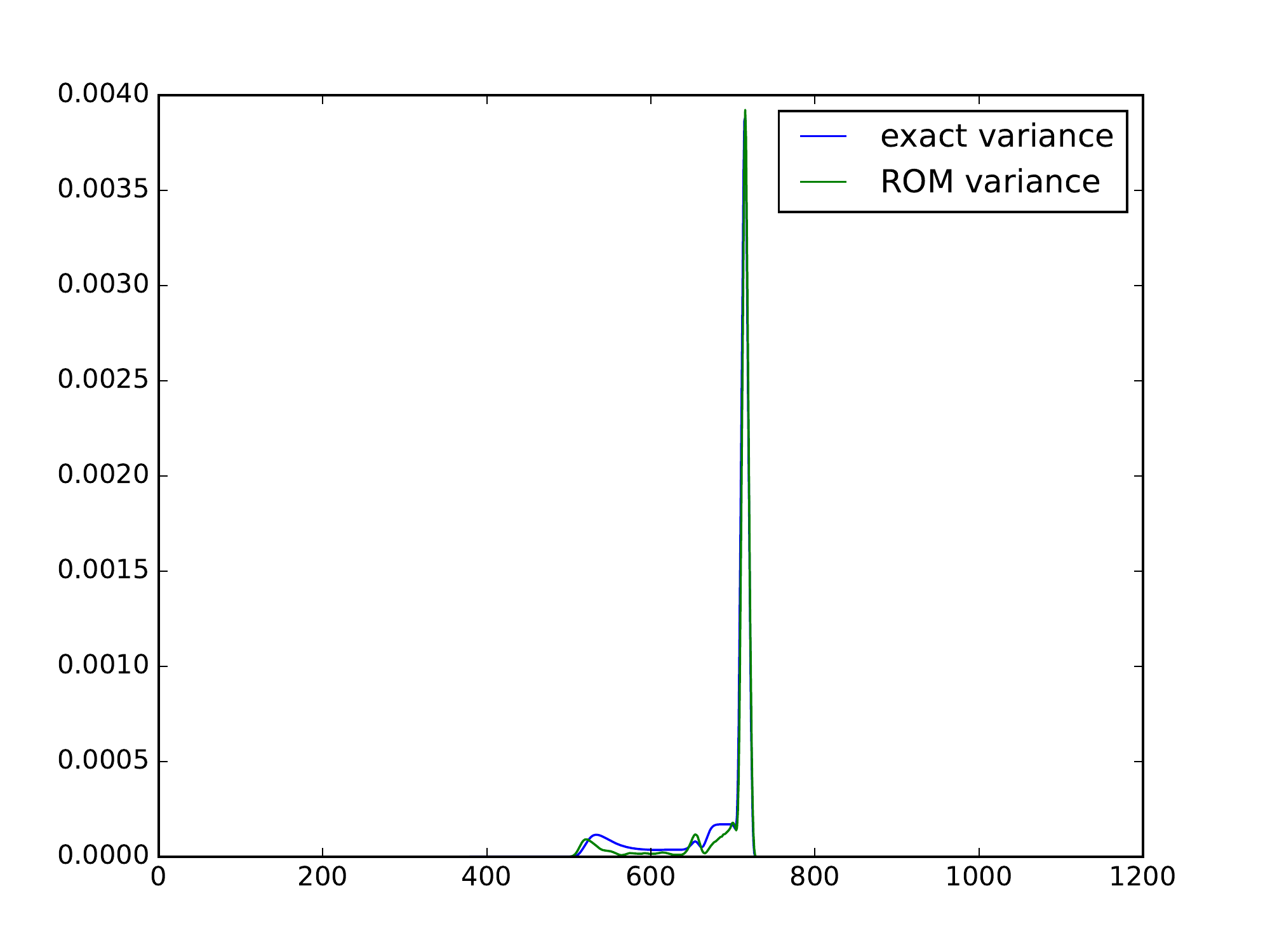}}
	\caption{\label{fig:mean_Euler_1d_2_param_mom} Solution mean, the mean plus/minus the standard deviation and the variance for both the reduced and the high-fidelity problem in the case of Euler equation with random initial condition and random flux for momentum}
\end{figure}

\begin{figure}[h!]
	\centering
	\subfigure[]{\includegraphics[width=0.49\textwidth,clip=]{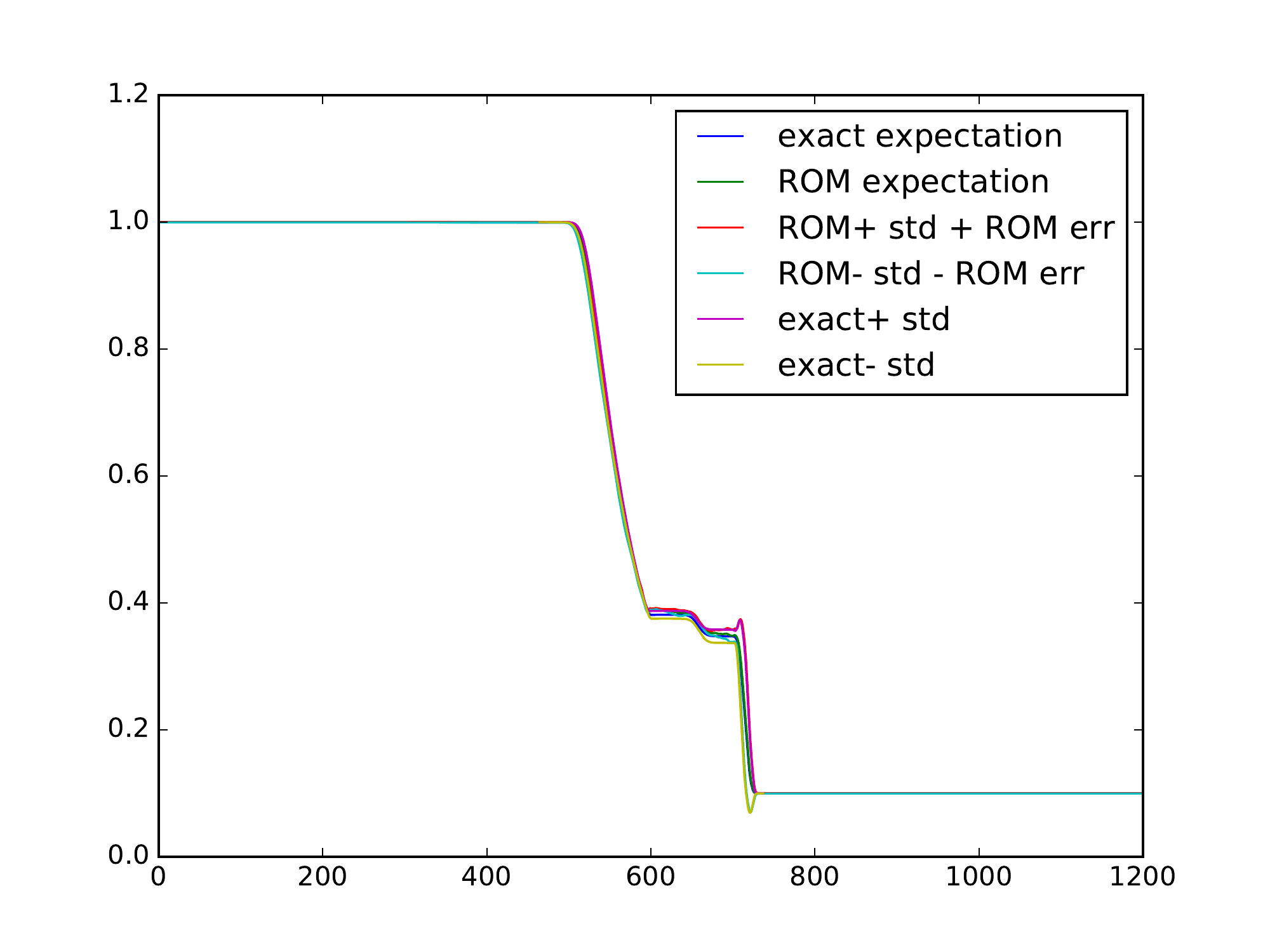}}
	\subfigure[]{\includegraphics[width=0.49\textwidth,clip=]{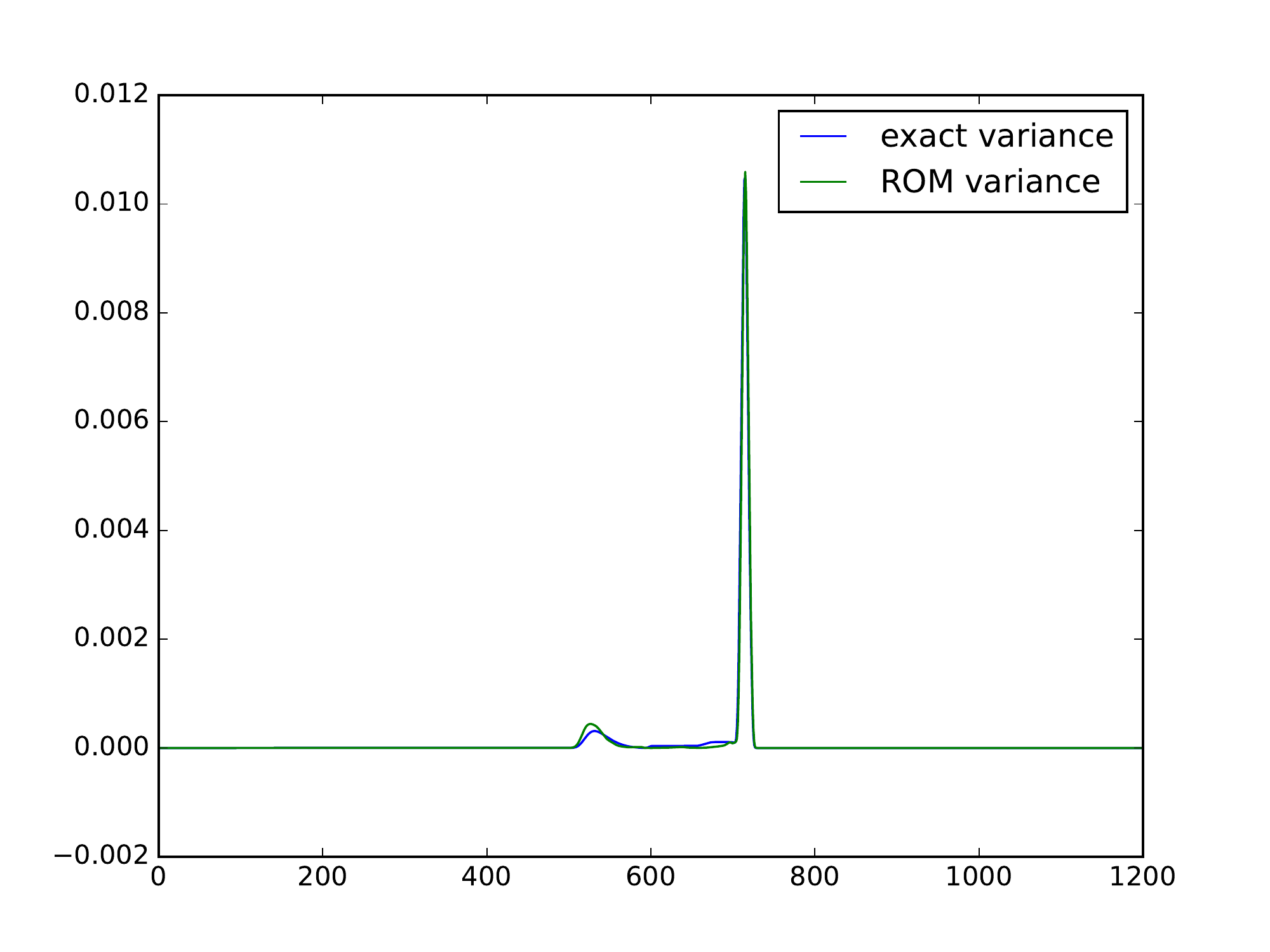}}
	\caption{\label{fig:mean_Euler_1d_2_param_E} Solution mean, the mean plus/minus the standard deviation and the variance for both the reduced and the high-fidelity problem in the case of Euler equation with random initial condition and random flux for the total energy}
\end{figure}

In the online phase, the UQ analysis is performed using a set with $100$ elements in the parameter domain $D_y=[-0.02,0.02] \times [1.4,1.5]$, which were generated by a random Monte Carlo method. Comparing again the solution mean and the variance, as well as the solution mean plus/minus the standard deviation of a random variable $\utruth{K}(w)$ in the case of the reduced problem and the high fidelity one (see Figures \ref{fig:mean_Euler_1d_2_param_rho}, \ref{fig:mean_Euler_1d_2_param_mom}, \ref{fig:mean_Euler_1d_2_param_E}), we obtain a computational saving time of 69\%. For a better visualization, we plot each component of the solution independently.  Indeed, the average computational time for one high fidelity solution is of $39.448$ seconds, while the reduced solution takes only $12.420$ seconds.

\subsection{Stochastic Sod's shock problem in 2D with random initial data and random flux}
Consider the two-dimensional Euler equations with random initial data and random flux:
\begin{align}\label{eq:Euler_2D}
	\frac{\partial \ubold}{\partial t}+\frac{\partial \fbold(\ubold,w)}{\partial x_1}+\frac{\partial \gbold(\ubold,w)}{\partial x_2}&=0, ~\xbold=(x_1,x_2) \in D=\lbrace (x_1,x_2) | x_1^2+x_2^2\leq 1\rbrace \\
	\ubold_0(\xbold,w)&=\ubold_0(\xbold,Y_1(w))
\end{align} 
where $y_j=Y_j(w), ~j=1,2, ~w \in \Omega$, the components are expressed as
$$\ubold=(\rho,\rho u,\rho v, E )^{T}, ~ \fbold=(\rho,\rho u^2+p,\rho u v,\rho u(E+p))^T, ~ \gbold=(\rho,\rho uv,\rho v^2+p,\rho v(E+p))^T$$ and the pressure as
$$p=(\gamma-1)\Big( E-\frac{1}{2}\rho (u^2+v^2)\Big).$$

We assume again randomness in the adiabatic constant, $\gamma=Y_2(w)$, and therefore $$\fbold(\ubold,w)=\fbold(\ubold,Y_2(w))$$ and $$\gbold(\ubold,w)=\gbold(\ubold,Y_2(w)).$$ The initial data is set in primitive variables as
\begin{equation*}
	\wbold_0(\xbold,w)=(\rho_0(\xbold,w),u_0(\xbold,w),v_0(\xbold,w),p_0(\xbold,w))^T=\begin{cases}
		(1,0,0,1), & \mbox{if } 0\leq r<0.5 \\ 
		(0.125+Y_1(w),0,0,0.1), & \mbox{if } 0.5<r\leq 1
	\end{cases}
\end{equation*}
where $r=\sqrt{x_1^2+x_2^2}$ is the distance of the point $(x_1,x_2)$ from the origin.

The computations have been performed on a triangular mesh consisting of 13548 cells and $N_h=6775$ DoFs, using $K=500$ time instances of step $\Delta t=0.0005$, the final time is $T=0.25$ and using a first order version of the RD scheme presented in \cite{abgrallLarat}. 

\begin{figure}[h!]
	\centering
	\includegraphics[width=0.52\textwidth,clip=]{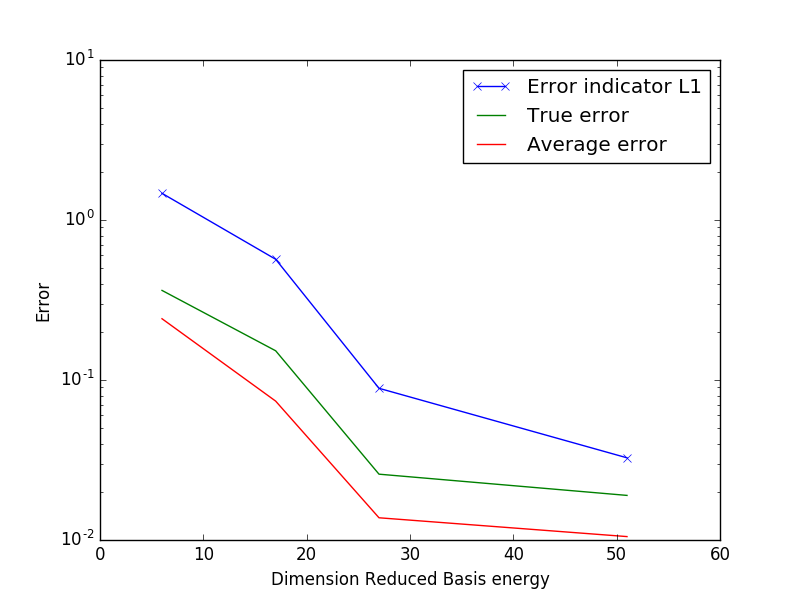}
	\caption{\label{fig:error_decay_2D} Error decay in Offline phase with respect to dimension of reduced basis space of Energy}
\end{figure}

In the offline step, the tolerance set for the greedy algorithm is $0.02$ and we are using a PODEIM--Greedy algorithm generating an EIM space with $(67, 68, 69, 76)$ basis functions and a RB space of dimension $(36, 50, 51, 53)$ in each component, namely in density, momentum in $x$ and $y$ direction and total energy. In this test case, we have randomness in both flux and initial condition, namely $Y_2(w)$, respectively $Y_1(w)$. We construct the random variables $Y_1(w),Y_2(w)$ using a uniform random Monte Carlo sampling method in the interval $D_y=[0.125,0.225] \times [1.4,1.6] $, resulting in a set with 100 elements. We can see the decay of the error during the Offline phase in \ref{fig:error_decay_2D}.

\begin{figure}[h!]
	\centering
	\subfigure[]{\includegraphics[width=0.49\textwidth,clip=]{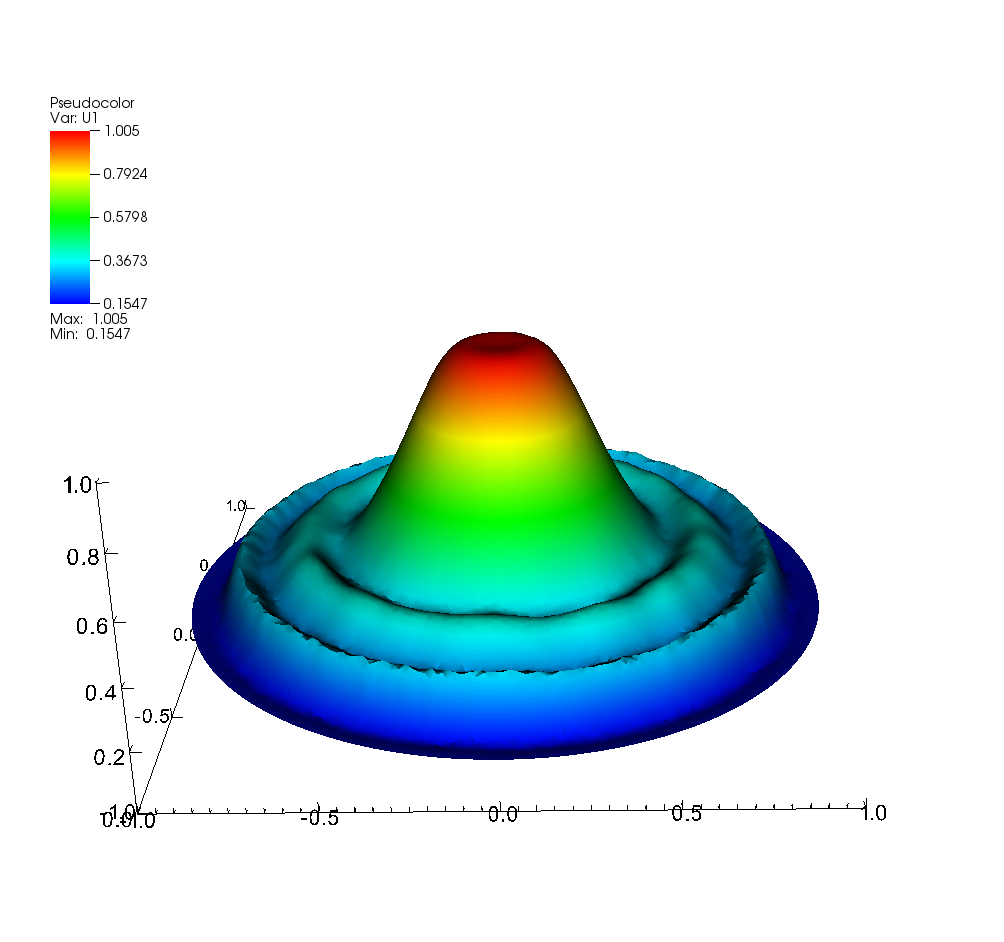}}
	\subfigure[]{\includegraphics[width=0.49\textwidth,clip=]{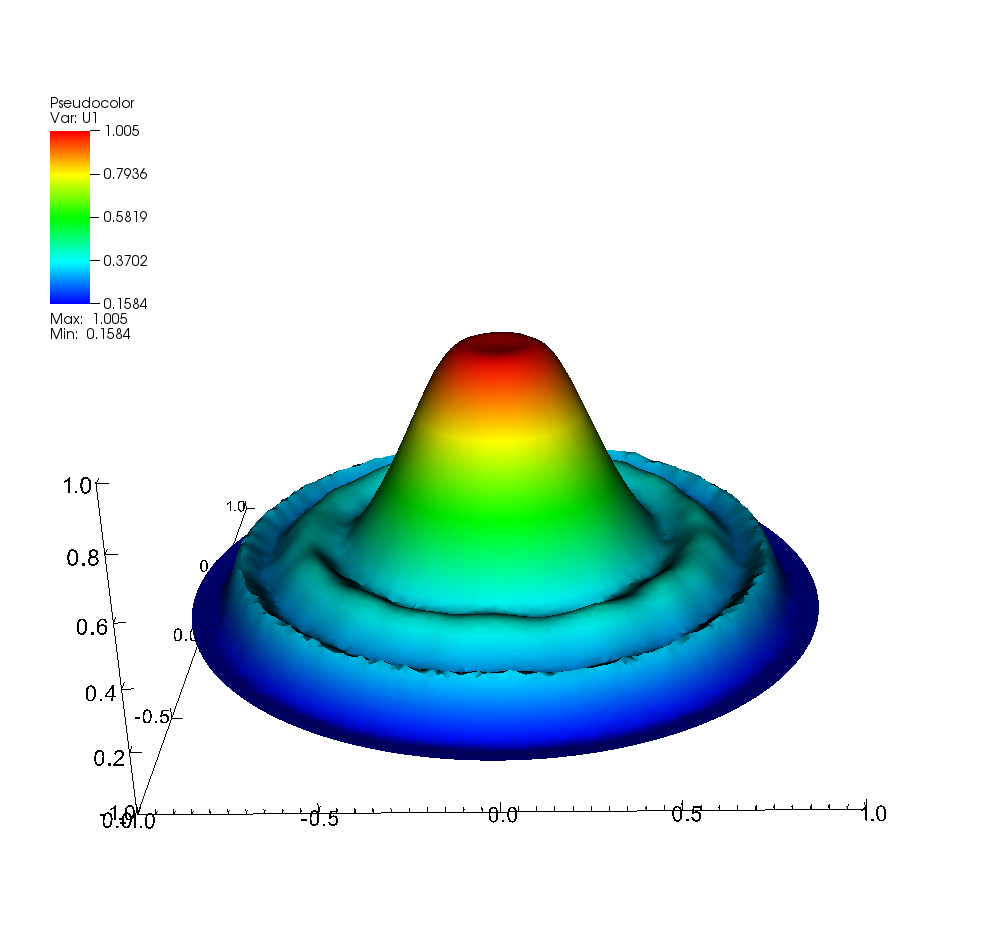}}
	\caption{\label{fig:RDRB} Density of high-fidelity solution (left) and the reduced solution (right) at final time T=0.25 for $Y=(0.16353811,  1.50632869)$ }
\end{figure}

\begin{figure}[h!]
	\begin{center}
		{\includegraphics[width=0.5\textwidth,clip=]{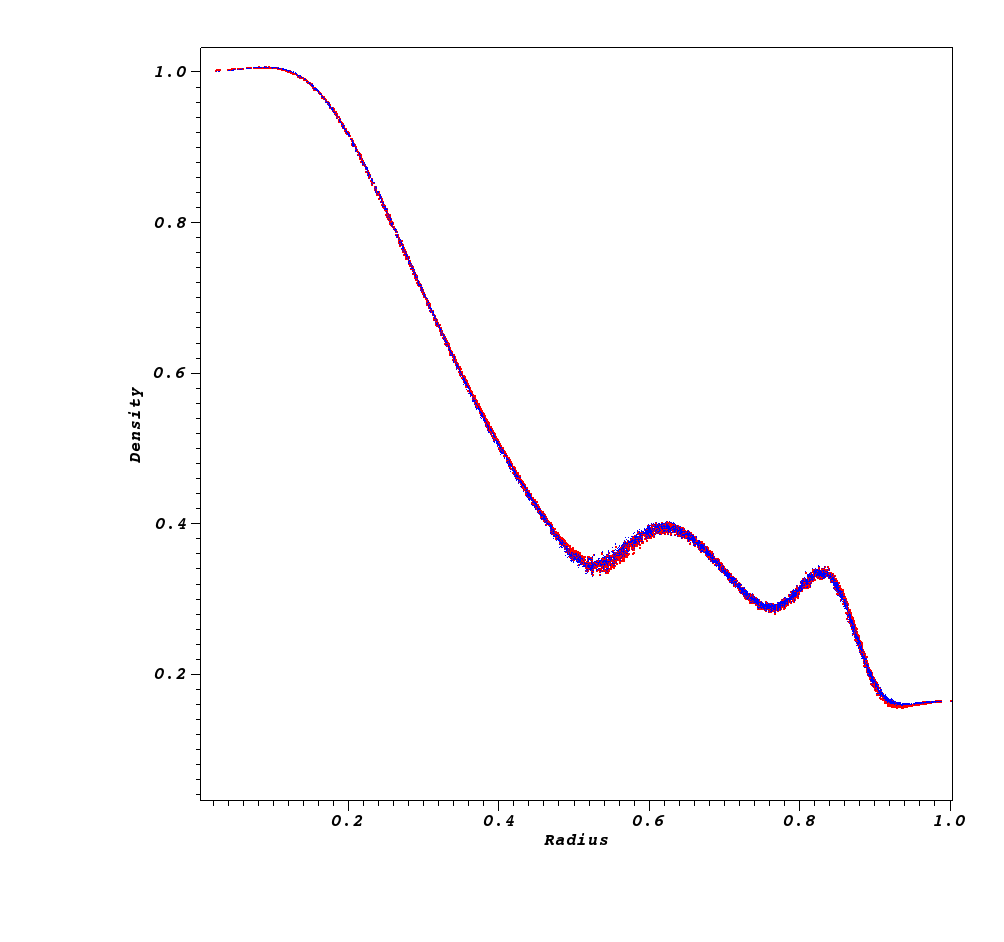}}
	\end{center}
	\caption{\label{fig:scatterRDRB} Scatter plot of density of the high-fidelity solution (red) and the reduced solution (blue) at final time T=0.25 for $Y=(0.16353811,  1.50632869)$  }
\end{figure}

\begin{figure}[h!]
	\centering
	\subfigure[]{\includegraphics[width=0.49\textwidth,clip=]{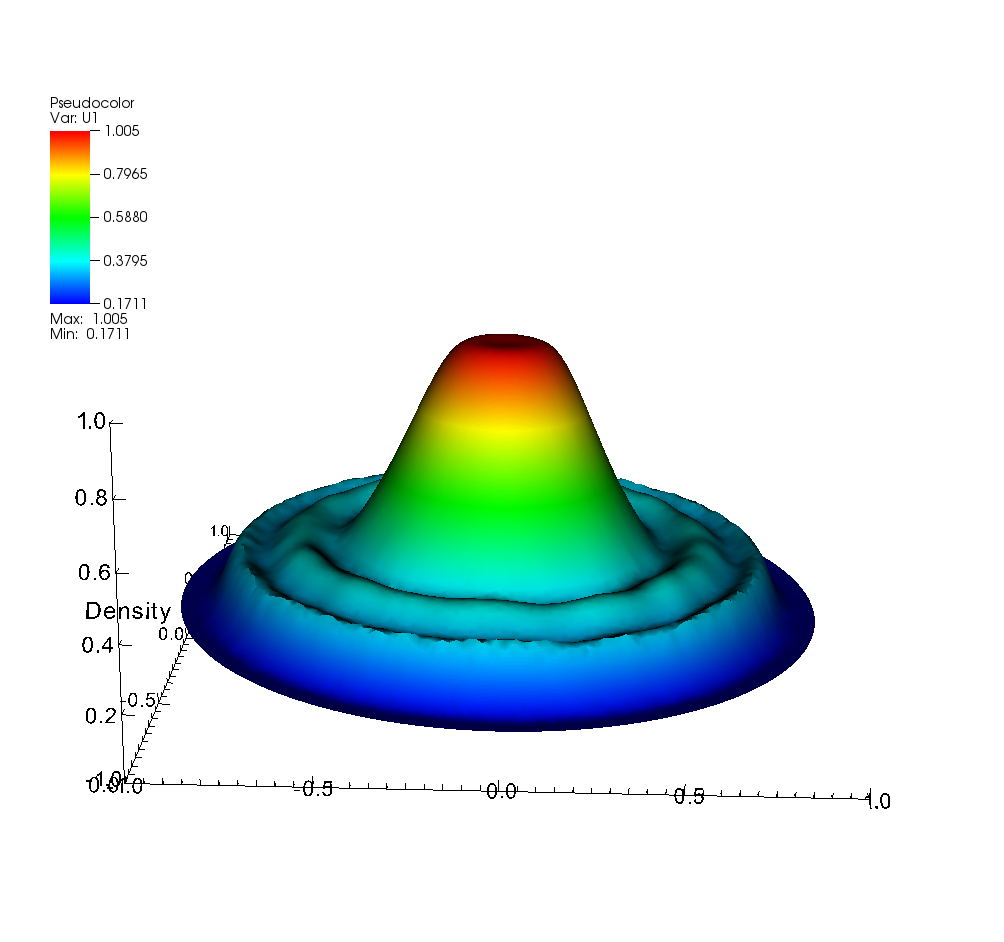}}
	\subfigure[]{\includegraphics[width=0.49\textwidth,clip=]{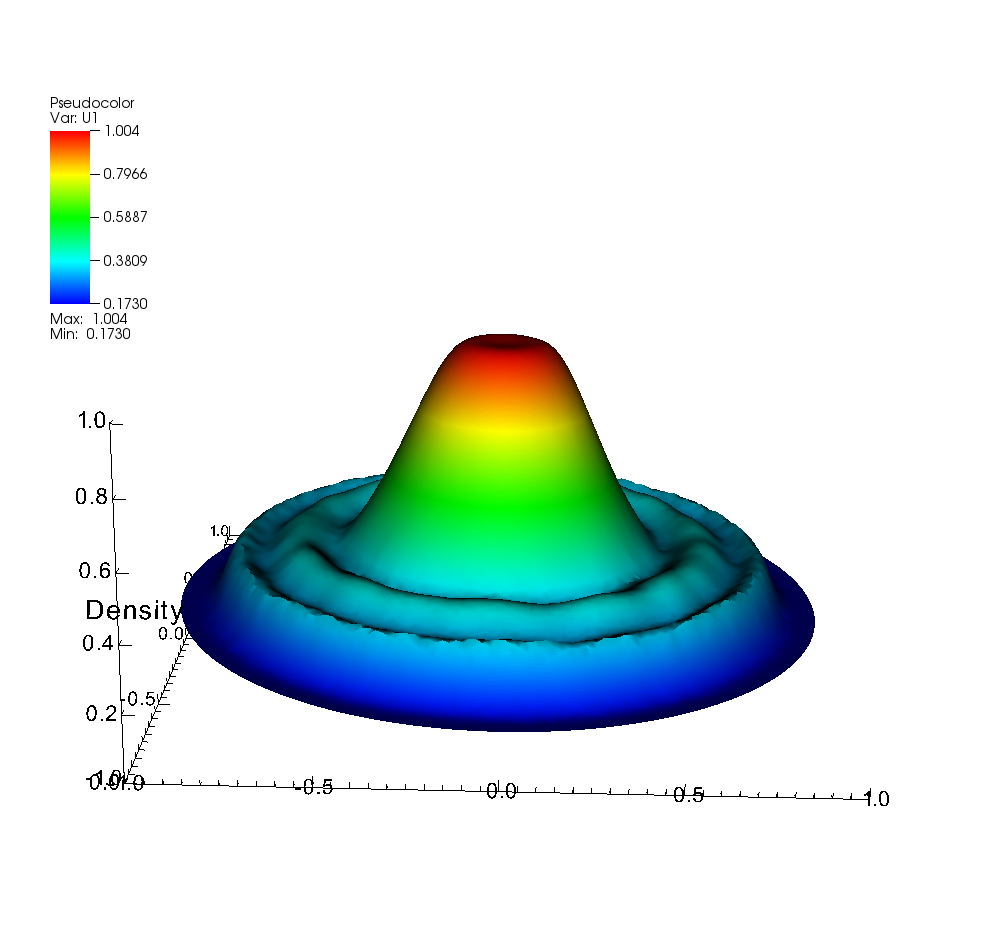}}
	\caption{\label{fig:meanRDRB} Solution mean for density of the high-fidelity problem (left) and  for the reduced solution (right) at final time T=0.25}
\end{figure}

\begin{figure}[h!]
	\begin{center}
		{\includegraphics[width=0.5\textwidth,clip=]{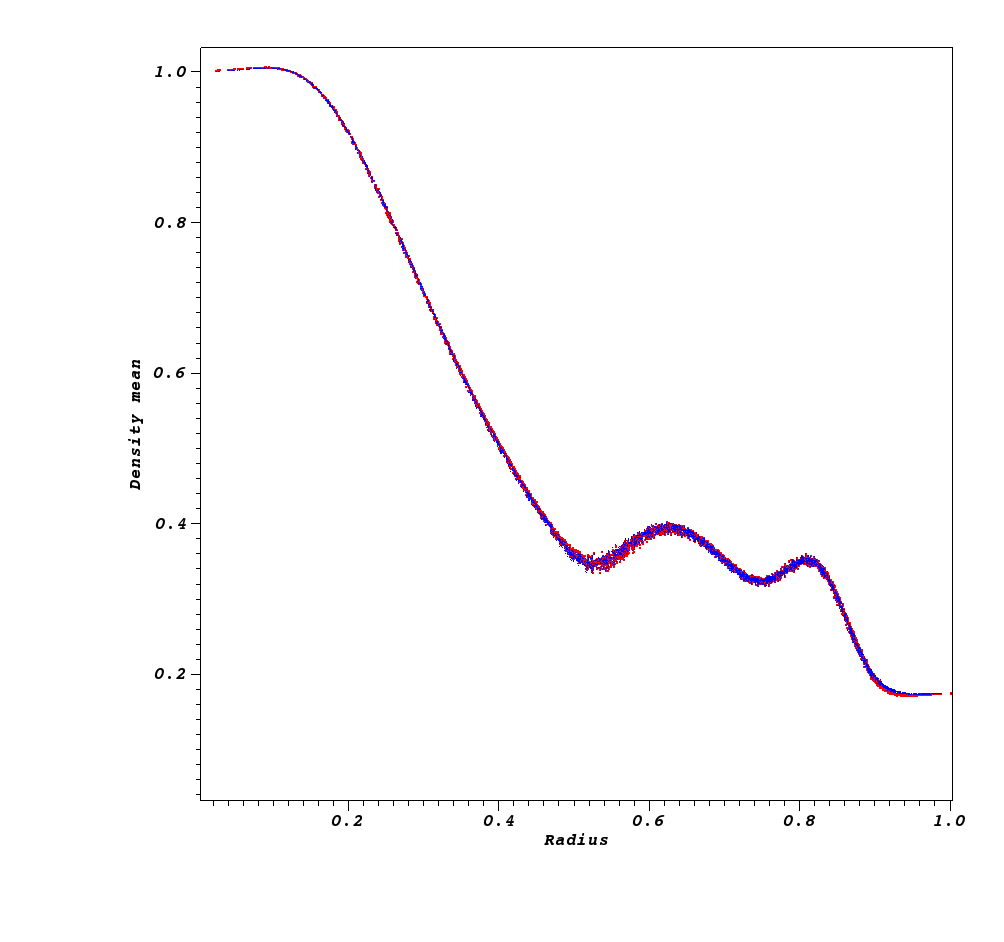}}
	\end{center}
	\caption{\label{fig:scattermeanRDRB} Scatter plot of density of the high-fidelity mean solution (red) and the mean of the reduced solution (blue) at final time T=0.25  }
\end{figure}

\begin{figure}[h!]
	\centering
	\subfigure[]{\includegraphics[width=0.49\textwidth,clip=]{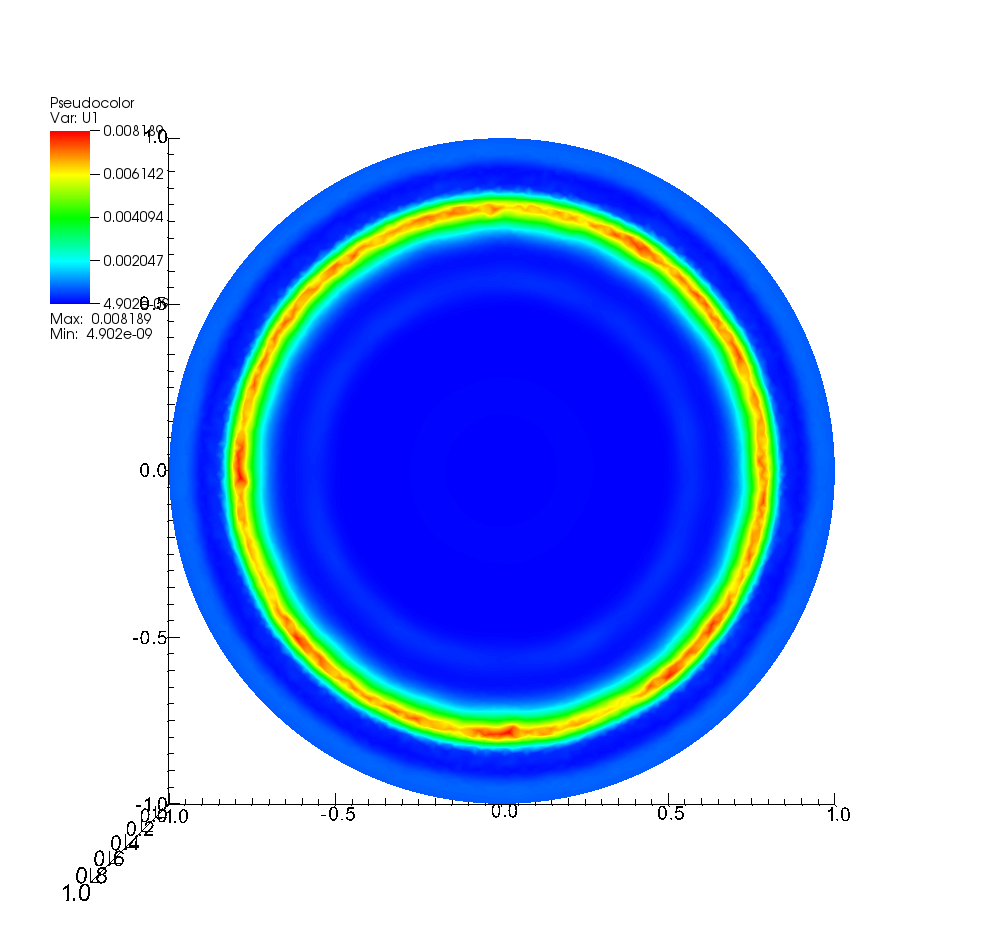}}
	\subfigure[]{\includegraphics[width=0.49\textwidth,clip=]{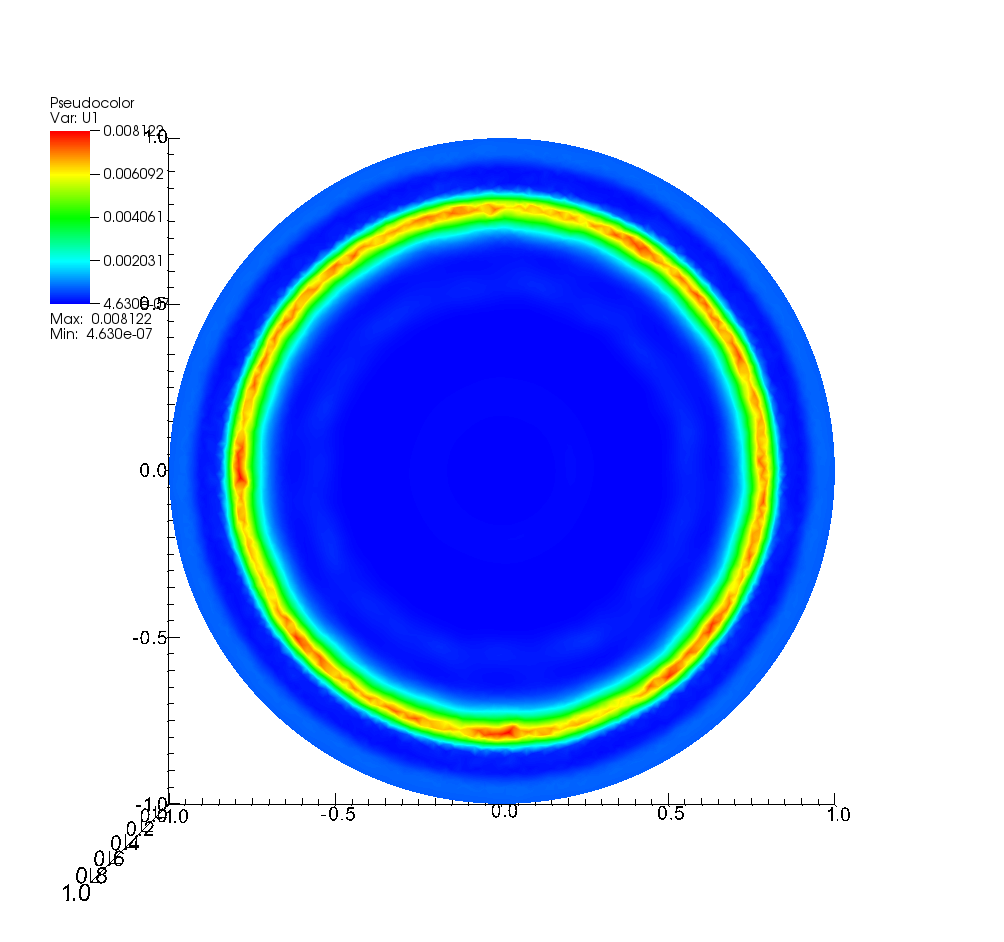}}
	\caption{\label{fig:varianceRDRB} Variance for the density of high-fidelity problem (left) and for the reduced solution (right) at final time T=0.25}
\end{figure}

\begin{figure}[h!]
	\begin{center}
		{\includegraphics[width=0.5\textwidth,clip=]{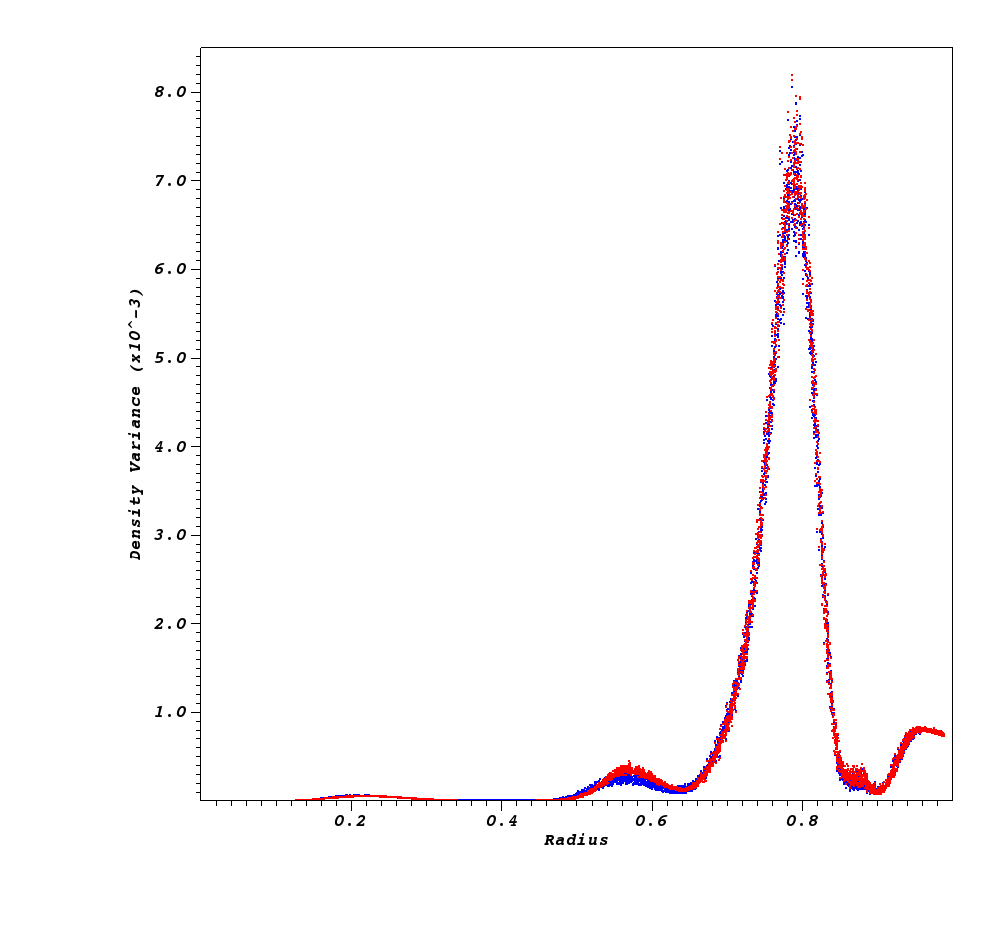}}
	\end{center}
	\caption{\label{fig:scattervarianceRDRB} Scatter plot of density of the high-fidelity variance (red) and the reduced solution variance (blue) at final time T=0.25  }
\end{figure}

In the online phase, the UQ analysis is performed using a set with 50 elements in the parameter domain $D_y=[0.125,0.225] \times [1.4,1.6]$, which were generated by a uniform random Monte Carlo method. Comparing again the solution mean (see Figures \ref{fig:meanRDRB}, \ref{fig:scattermeanRDRB}) and the variance (see Figure \ref{fig:varianceRDRB}, \ref{fig:scattervarianceRDRB}), in the case of the reduced problem and the high fidelity one (see Figure \ref{fig:RDRB}, \ref{fig:scatterRDRB}), we can see that the reduced solution has qualitatively no worsening. Morover, we obtain a computational saving time of 76\%. Indeed, the average computational time for one high fidelity solution is of 517.59 seconds, while the reduced solution takes only 125.50 seconds.

\section*{Acknowledgements}
In this work, R.A and S. T. have been funded in part by the SNF project 200021\_153604 "High fidelity simulation for compressible materials".  R.C. has been funded by the University of Z\"urich. D.T. has been funded by the ITN project "ModCompShock: Computational modelling of shocks and interfaces" funded  SERI agreement SBFI Nr 15.0269-1.

\newpage

%% file: Bibliography.tex
\bibliographystyle{siam}
\bibliography{biblio,biblio_ch_uq,biblio_UQ,biblio_algo}{}